\documentclass[12pt]{amsart}

\usepackage{amssymb}
\usepackage{amsmath}
\evensidemargin\oddsidemargin

\begin{document}

\setcounter{page}{1}

\newtheorem{theorem}{Theorem}

\newtheorem{THEO}{Theorem\!\!}
\newtheorem{Conj}{Conjecture}

\renewcommand{\theTHEO}{}

\newcommand{\eqnsection}{
\renewcommand{\theequation}{\thesection.\arabic{equation}}
    \makeatletter
    \csname  @addtoreset\endcsname{equation}{section}
    \makeatother}
\eqnsection

\def\a{\alpha}
\def\b{\beta}
\def\Ab{A^\b}
\def\CC{{\mathbb{C}}} 
\def\cD{{\mathcal{D}}} 
\def\cH{{\mathcal{H}}} 
\def\cL{{\mathcal{L}}} 
\def\cN{{\mathcal{N}}} 
\def\cZ{{\mathcal{Z}}} 
\def\Ea{E_\a}
\def\Eac{{\mathcal E}}
\def\EE{{\mathbb{E}}} 
\def\elaw{\stackrel{d}{=}}
\def\eps{\varepsilon}
\def\Fa{F_\a}
\def\Ga{G_\a}
\def\hT{{\hat T}}
\def\hX{{\hat X}}
\def\ii{{\rm i}}
\def\lbd{\lambda}
\def\lacc{\left\{}
\def\lcr{\left[}
\def\lpa{\left(}
\def\lva{\left|}
\def\NN{{\mathbb{N}}} 
\def\pb{{\mathbb{P}}}
\def\rl{{\mathbb{R}}}
\def\racc{\right\}}
\def\rcr{\right]}
\def\rpa{\right)}
\def\rva{\right|}
\def\Ta{T_\a}
\def\Ua{U_\a}
\def\Un{{\bf 1}}
\def\V{{\rm Var}}
\def\ZZ{{\mathbb{Z}}} 

\def\d{\, \mathrm{d}}
\def\qed{\hfill$\square$}
\def\elaw{\stackrel{d}{=}}
\def\claw{\stackrel{d}{\rightarrow}}

\newcommand{\fin}{\vspace{-0.3cm}
                  \begin{flushright}
                  \mbox{$\Box$}
                  \end{flushright}
                  \noindent}

\title[Persistence probabilities \& exponents]{Persistence probabilities \& exponents}

\author[F. Aurzada]{Frank Aurzada}
\address{Technische Universit\"at Berlin, Institut f\"ur Mathematik, Strasse des 17. Juni 136, D-10623 Berlin. {\em Email}: {\tt aurzada@math.tu-berlin.de}}
\author[T. Simon]{Thomas Simon}
 \address{Laboratoire Paul Painlev\'e, Universit\'e Lille 1, Cit\'e Scientifique, F-59655 Villeneuve d'Ascq Cedex. {\em Email}: {\tt simon@math.univ-lille1.fr}}

\begin{abstract}
This article deals with the asymptotic behaviour as $t\to +\infty$ of the survival function $\pb[T > t],$ where $T$ is the first passage time above a non negative level of a random process starting from zero. In many cases of physical significance, the behaviour is of the type $\pb[T > t]=t^{-\theta + o(1)}$ for a known or unknown positive parameter $\theta$ which is called a persistence exponent. The problem is well understood for random walks or L\'evy processes but becomes more difficult for integrals of such processes, which are more related to physics. We survey recent results and open problems in this field.        
\end{abstract}



\maketitle

\tiny






\normalsize

\tableofcontents


\section{Introduction}

Let $\{X_t, \, t\ge 0\}$ be a real stochastic process in discrete or continuous time, starting from zero. The analysis of the first passage time $T_x = \inf \{t > 0, \, X_t > x\}$ above a non-negative level $x$ is a classical issue in probability. In this paper we will be concerned with the asymptotic behaviour as $t\to +\infty$ of the survival function $\pb[T_x > t]$ for a class of processes related to random walks and L\'evy processes. This problem has attracted some interest in the recent literature under the denomination persistence probability. In the self-similar framework, it is also related to the lower tail probability problem which is the asymptotic study as $\varepsilon\to 0$ of the quantity $\pb [S_t <\varepsilon]$ where $S_t = \sup\{X_s, s\le t\}$ is the one-sided supremum process. In many situations of interest it turns out that the behaviour is polynomial: one has
\begin{equation}
\label{PolBev}
\pb [T_x > t]\; = \; t^{-\theta +o(1)}
\end{equation}
for a non-negative parameter $\theta$ called a persistence exponent, which usually does not depend on $x >0$ and often belongs to $[0,1].$ The study of asymptotic behaviours like (\ref{PolBev}) has gained some attraction over the last years in the physical literature as well, where the parameter $\theta$ is often called a survival exponent. Notice that estimates of the type (\ref{PolBev}) also appear in reliability theory, a subject that we shall here not discuss, where $T_x$ is viewed as a certain failure time whose typical upper tails are Pareto-like - see \cite{KP}.
 
If $X$ is a random walk or a L\'evy process, studying the law of $T_x$ is a special part of fluctuation theory and the persistence probabilities are then well understood. For example if $X_1$ is attracted to a stable law, classical fluctuation identities entail that the persistence exponent is the positivity parameter of the latter. There are many accounts on fluctuation theory and we refer in particular to \cite{Bi3} for random walks and to \cite{Do} for L\'evy processes. In the first part of this article we focus on the results of this theory dealing with the asymptotic behaviour of $\pb[T_x > t],$ and we try to be as exhaustive as possible. In general, the behaviour is the same in discrete and continuous time, except of course for $x = 0$ where the problem becomes different for L\'evy processes. We notice that even though a posteriori the resulting exponents turn out to be the same, there is no simple approximation argument that would yield this a priori, and for this reason we have to consider the discrete and continuous time situations separately. The results, all classical, are presented here in order to give some insight into more complex situations where $X$, though constructed upon a L\'evy process or a random walk, is not a Markov process anymore.

Some of these more complex situations, which are called non-trivial in the physical literature, are the matter of the second part of this article. We first consider integrated random walks and integrated L\'evy processes. It turns out that for such simple constructions, the computation of the persistence exponent is not quite easy in general. We display recent and less recent results where the persistence probabilities are estimated with various degrees of precision, and which all have the common feature that the persistence exponent of the underlying process is twice the persistence exponent of its integral. We believe that this is a kind of universal rule and several conjectures are stated in this direction, the most general being Conjecture 3 and Conjecture 5. A curious fact is that the situation with no positive jumps, where classical fluctuation theory becomes much simpler, appears to be more difficult than the dual situation with no negative jumps, where some exponents of integrated processes can be computed. The reason for that seems to be only methodological, and to cast out the curse of spectral positivity in the study of persistence probabilities of integrated processes probably requires other, say deeper, tools. 

We then consider fractionally integrated processes, where the situation is more difficult. The case when the underlying L\'evy process is Brownian motion yields self-similar Gaussian processes which can be changed into a stationary one by the Lamperti transformation. The sought after persistence exponent is then directly related to an estimate on the probability of non-zero crossings for these Lamperti transforms, a problem known to be hard in spite of the It\^o-Rice formula which gives some information at the expectation level in the smooth case. We present some universality and monotonicity results, and also some partially heuristic comparisons with the fractional Brownian motion case, whose persistence exponent can be computed explicitly. We also display other explicit computations of persistence exponents for related processes such as weighted sums or iterated processes.

The third and last part of this article deals with some applications of the persistence probabilities in mathematical physics. We first deal with Lagrangian regular points of the inviscid Burgers equation with random self-similar data. The link between the Hausdorff dimension of such points and the persistence exponent of integrated processes dates back to the original paper \cite{Si1}. We recall here that the problem is still open in the fractional Brownian motion case and state a plausible conjecture when the initial data is a two-sided stable L\'evy process. Second, we consider the zero-crossings of a peculiar Gaussian stationary process which is related to the positivity of Kac polynomials with large even degree. This connection, which was discovered in \cite{DPSZ}, has further ramifications with the persistence exponents of integrals of Brownian motion with higher order and of a certain diffusion equation with white noise initial conditions in the plane, and we make a brief account on the subject. We last consider three different interacting statistical systems whose analysis hinges significantly upon the persistence probabilities of integrated processes: wetting models with Laplacian interaction, fluctuating interfaces with Langevin dynamics, and sticky particles on the line with Poissonian initial conditions. 

Some open problems stated in the present paper are believed to be challenging and we think that they could catch the attention of some colleagues. We finally point out that we have not exhausted here all implications of persistence in physics and that the persistence exponent of many other models remains unknown - see \cite{CDPMS, M} and the references therein. 

\section{Classical results}

\subsection{Random walks} 

\label{sec:rw}

Let $\lacc X_n, \, n\ge 1\racc$ be a sequence of i.i.d. real random variables with common distribution $\mu$, and $S_n = X_1 + \ldots + X_n$ be the associated random walk. Consider
$$T_x \; =\; \inf\lacc n\ge 1, \; S_n > x\racc$$
the first-passage time above $x \ge 0.$ Recall the following basic rule, to be found in e.g. Chapter XII.1-2 of \cite{F}:
\begin{equation}
\label{Trans}
T_x \; \mbox{is a.s. finite for every $x\ge 0$}\;\Leftrightarrow\; T_0 \; \mbox{is a.s. finite},\end{equation} 
and that the latter is also equivalent to the fact that $S_n$ does not drift to $-\infty.$ In this case one has $\pb[T_x > n]\to 0$ as $n\to +\infty$ for every $x\ge 0,$ and the difficulty to estimate the speed of convergence comes from the fact that the event $\{T_x > n\} = \{S_1 \le x, \ldots, S_n \le x\}$ depends on $n$ correlated random variables.  \\


If $\mu$ is concentrated on $\rl^+,$ then $\pb[T_x > n] = \pb[S_n < x]$ and the problem becomes one-dimensional. The straightforward inequality
$$\pb[S_n < x] \; \le\; e^{x + n\log (\EE[e^{-X_1}])}$$
shows that $\pb[T_x > n]$ tends to zero at least exponentially fast (unless $\mu$ is degenerate at zero).  Notice that the rate might also be superexponential and depend on $x$ at the logarithmic scale: if $X_1$ has a positive strictly $\a-$stable law for instance, de Bruijn's Tauberian theorem - see Theorem 4.12.9 in \cite{BGT} - leads to 
$$-\log \pb[T_x > n]\; =\; -\log \pb[S_n < x]\; =\; -\log \pb[S_1 < x n^{-1/\a}]\;\sim\; \kappa_\a x^{-\a/(1-\a)} n^{1/(1-\a)}$$
for every $x > 0,$ with some explicit $\kappa_\a >0.$ Of course, one has $T_0 = 1$ a.s. whenever $\mu$ does not charge 0. If $\mu\{0\} > 0$, Jain and Pruitt's general uniform results on renewal sequences express the asymptotic behavior of $\pb[S_n < x]$ in terms of quantities related to $\mu$, at the logarithmic scale  - see Theorem 2.1 in \cite{JP} - and also at the exact scale, under some conservativeness assumption  - see Theorem 4.1 in \cite{JP}. \\

If $\mu$ is not concentrated on $\rl^+,$ it is easy to see that the asymptotics of $\pb[T_x > n]$ will not drastically depend on $x$: choosing $\eps > 0$ such that $\mu_\eps=\mu(-\infty, -\eps) > 0,$ the Markov property entails $\pb[T_x > n] \ge \pb[T_0 > n] \ge  \pb[T_\eps > n-1]\pb[X_1<-\eps]\ge \mu_\eps\pb[T_\eps > n]\ge \mu_\eps^k\pb[T_x > n]$ as soon as $k\eps \ge x,$ so that 
$$\pb[T_x > n] \;\asymp\; \pb[T_0 > n]$$ 
for every $x \ge 0.$ \\

The following formula computes the generating function of $\{ \pb[T_0 > n], \; n\ge 0\}$ in terms of that of the sequence $\{n^{-1}\pb [S_n \le 0], \; n\ge 1\},$ and is true for any random walk. 

\vspace{2mm}

\noindent
{\bf Sparre Andersen's formula.} For every $z \in ]-1,1[$ one has
\begin{equation}
\label{SA}
\sum_{n\ge 0} z^n \pb[T_0 > n]\; = \; \exp\lcr \sum_{n\ge 1} \frac{z^n}{n} \pb [S_n \le 0]\rcr.
\end{equation}

\vspace{2mm}
This result is obtained after simple rearrangements from Theorem XII.7.1 in \cite{F}, whose proof has a combinatorial character and depends heavily on the independence and stationarity of the increments of the random walk. A simpler method relying on elementary Fourier analysis - see Chapter 3.7 in \cite{DMcK} - yields the more general

\vspace{2mm}

\noindent
{\bf Spitzer's formula.} For every $z \in ]-1,1[$ and $\lbd \ge 0$, one has
\begin{equation}
\label{SP}
\sum_{n\ge 0} z^n \EE[e^{-\lbd M_n}]\; = \; \exp\lcr \sum_{n\ge 1} \frac{z^n}{n} \EE [e^{-\lbd S_n^+}]\rcr
\end{equation}
with the notation $M_n = \max (0, S_1, \ldots, S_n)$ and $S_n^+ = \max (0, S_n)$.
\vspace{2mm}

Indeed, one obtains (\ref{SA}) as a consequence of (\ref{SP}) in letting $\lbd \to +\infty.$ It is beyond the peculiar scope of the present paper to discuss the full strength and the various generalisations of Spitzer's formula such as Baxter-Spitzer's formula or the Wiener-Hopf factorization, and we refer e.g. to \cite{Bi3, Do, F} for more on this subject. Let us simply remark that Sparre Andersen's formula entails easily - see e.g. Theorem XII.7.2 in \cite{F} - the following characterization of (\ref{Trans}):
$$S_n \; \mbox{drifts to $-\infty$}\;\Leftrightarrow\; \sum_{n\ge 1} \frac{1}{n}\, \pb[S_n > 0]\; <\; +\infty.$$ 
From now on, we will suppose that $S_n$ does not drift to $-\infty$ and that $\mu$ is not concentrated on $\rl^+$. This entails that $\pb[S_n> 0] >0$ and $\pb[S_n< 0] >0$ for every $n\ge 1$.  \\

A remarkable consequence of (\ref{SA}) is that when $\pb [S_n \le 0] =\rho\in (0,1)$ for every $n\ge 1,$ one obtains an explicit formula for $\pb[T_0 > n]$ depending only on $\rho$ and $n:$ one has
$$\sum_{n\ge 0} z^n \pb[T_0 > n]\; = \; \frac{1}{(1-z)^\rho}\; =\; \sum_{n\ge 0} \frac{\Gamma(n+\rho)}{n!\Gamma(\rho)}\, z^n,$$
so that
\begin{equation}
\label{Exact}
\pb[T_0 > n]\; = \; \frac{\Gamma(n+\rho)}{n!\Gamma(\rho)}\; \sim\; \frac{n^{\rho-1}}{\Gamma(\rho)}\cdot
\end{equation}
For example, symmetric random walks such that $\pb [S_n = 0] = 0$ for every $n\ge 1$ (this latter property is true when $\mu$ is non atomic, for instance) all enjoy the property that
$$\pb[T_0 > n]\; = \; \frac{\Gamma(n+1/2)}{\sqrt{\pi} n!}\; \sim\; \frac{1}{\sqrt{\pi n}}\cdot$$
Recall in passing that the estimate is slightly different for simple random walks with $\mu\{1\} = \mu\{-1\} = 1/2$, since then $\pb[S_{2n}=0]\neq 0.$ The classical computation 
$$\sum_{n\ge 0} z^n \pb[T_0 > n] \; =\; \lpa \frac{1}{1-z} - \frac{1-\sqrt{1-z^2}}{z(1-z)}\rpa\;\sim\; \sqrt{\frac{2}{1-z}}\;\; \mbox{as $z\uparrow 1,$}$$
and the Tauberian theorem for monotonic sequences entail $\pb[T_0 > n]\; \sim\; \sqrt{2/\pi n}.$ \\

Another remarkable consequence of (\ref{SP}) is the exact computation of the persistence exponent whenever $\{S_n, n\ge 1\}$ fulfils the so-called Spitzer's condition
\begin{equation}
\label{Spitz}
\lim_{n\to \infty}\pb[S_n < 0]\; =\; \rho\in (0,1).
\end{equation}
 A theorem of Rogozin shows indeed that (\ref{Spitz}) entails
\begin{equation}
\label{Rogo}
\pb[T_x > n]\; \sim\;  \frac{c_x n^{\rho-1} l(n)}{\Gamma(\rho)} \; =\; n^{\rho -1 + o(1)}
\end{equation}
for every $x\ge 0$ with $l(n)$ some slowly varying sequence and $c_x$ some explicit positive constant. Besides, as explained in Theorem 8.9.12 of \cite{BGT}, this asymptotic behavior is actually equivalent to (\ref{Spitz}). Recall - see Chapter 7 in \cite{Do} - that  (\ref{Spitz}) is also equivalent to $\pb[S_n <0] \to \rho$ in Ces\`aro mean.
When $\pb [S_n < 0]\to 0$ various behaviors are possible, contrary to the above. For example if $n^\rho \pb [S_n < 0]$ is slowly varying for some $\rho \in (0,1),$ Theorem 8.9.14 in \cite{BGT} and Theorem XVII.5.1 in \cite{F} yield
\begin{equation}
\label{Embre}
\pb[T_0 > n]\; \sim\; l(n) n^{-(1+\rho)}
\end{equation}
with $l(n)$ some slowly varying sequence. Other various behaviours also appear when $\mu$ has positive expectation and we refer to \cite{Do89} for precise results. When $\pb [S_n < 0]\to 1,$ there does not seem to be any universal behavior for $\pb[T_0 > n]$, although one might expect that the persistence exponent is always zero. \\

We conclude this paragraph with the so-called upward skip-free or right-continuous random walks on $\ZZ$ viz. such that Supp $\mu\,\subset\, \{1, 0, -1, -2, -3, \ldots\}$, for which our survival analysis does not require the use of (\ref{SA}). Indeed, the distribution of $T_k$ is then given by Kemperman's formula \cite{K} which reads: 
\begin{equation}
\label{Kempe}
\pb [T_k = n] \; =\; \frac{k+1}{n}\, \pb[S_n = k+1]\quad \mbox{for every $n > k.$}
\end{equation}
In this particular case, we are reduced to the asymptotical behavior of the one-dimensional probability $\pb[S_n = 1].$ Notice - see \cite{Ke} - that there is a universal upper bound $\pb[S_n = 1]\le c n^{-1/2}$ leading to $\pb [T_k > n]\le 2c (k+1) n^{-1/2},$ but the exact behavior of $\pb[S_n = 1]$ depends on $\mu$. For example, if $\mu$ is in the respective domain of normal attraction of some strictly $\a-$stable law - with $\a\in(1,2]$ since otherwise $S_n$ would drift to $-\infty$ by the law of large numbers, Gnedenko's local limit theorem - see e.g. Theorem 8.4.1 in \cite{BGT} - yields $\pb[S_n = 1] \sim c_\a n^{-1/\a}l(n)$ for some explicit $c_\a$, so that by (\ref{Kempe})
\begin{equation}
\label{Gned}
\pb[T_k > n]\;\sim\; \a c_\a (k+1) l(n) n^{-1/\a}.
\end{equation}
Recall - see Proposition 8.9.16 in \cite{BGT} - that $\mu$ is in the domain of attraction of some strictly $\a-$stable law with $\a\in(1,2]$ if and only if (\ref{Spitz}) holds with $\rho = 1 - 1/\a,$ so that (\ref{Gned}) is actually a consequence of (\ref{Rogo}). On the other hand, the slowly varying term $l(n)$ can be removed if $S_n$ is in the respective domain of normal attraction (viz. when $S_n/n^{1/\a}$ converges in law to some non degenerate limit - see e.g. the concluding  remark of Chapter XVII.5 in \cite{F}), and this degree of precision is not given by Rogozin's theorem.

\subsection{L\'evy processes} 

\label{sec:lp}

Let $\lacc Z_t, \; t\ge 0\racc$ be a non-degenerate real L\'evy process starting from 0 and 
$$T_x \; =\; \inf\lacc t > 0, \; Z_t > x\racc$$
be its first passage time above $x\ge 0.$ Consider $\lacc Z_n, \; n\ge 1\racc$ the associated random walk. The inequality
\begin{equation}
\label{Rough}
\pb [T_x > t]\; \le \; \pb[{\tilde T}_x > [t]]
\end{equation}
with the notation ${\tilde T}_x = \inf\lacc n \ge 1, \; Z_n > x\racc,$ yields a rough upper bound for $\pb [T_x > t]$ which can be made more precise as a function of $t$ and $x$  in applying the results of the previous paragraph. This upper bound however does not yield enough information in general if $x =0.$ For example Rogozin's criterion - see Proposition VI.11 in \cite{Be} - shows that $T_0 = 0$ a.s. when $Z$ has unbounded variations, whereas the function $t\mapsto\pb[Z_1 \le 0, \ldots, Z_{[t]}\le 0]$ has a positive limit at $+\infty$ if $Z$ drifts to $-\infty.$ Recall also that in the bounded variation case, the regularity of the half-line for $Z$ is characterized in terms of the L\'evy measure - see Theorem 22 in \cite{Do}, and again $\pb[Z_1 \le 0, \ldots, Z_{[t]}\le 0]$ might have a positive limit when $t \to +\infty$ eventhough $T_0 = 0$ a.s. \\

If $x >0$ however, it turns out that the two quantities in (\ref{Rough}) are often comparable. First of all, for every $x >0$ one has 
\begin{equation}
\label{Defect}
\pb[T_x = +\infty]\, >\, 0\; \Leftrightarrow\; Z_t\to -\infty\;\mbox{a.s.}\; \Leftrightarrow\; \int_1^\infty \frac{1}{t} \,\pb [Z_t > 0] dt \, <\, +\infty
\end{equation}
(see e.g. Proposition 4.6 in \cite{Do}), and a simple analysis shows that the latter condition is equivalent to
$$Z_n\to -\infty\;\mbox{a.s.}\; \Leftrightarrow\;\sum_{n\ge 1}\frac{1}{n} \,\pb [Z_n > 0]\, <\, +\infty,$$ 
so that with the above notation $\pb[T_x > t] \to 0$ if and only if $\pb[{\tilde T}_x > t] \to 0$ as $t\to +\infty.$ In the following, we will see that in many examples one has 
$$\pb [T_x > t]\; \asymp \; \pb[{\tilde T}_x > [t]].$$
However, it does not seem easy to prove this estimate a priori, which probably does not hold in full generality. \\

We will suppose henceforth that (\ref{Defect}) does not hold. Set $\mu$ for the law of $Z_1$, $M_t = \sup\lacc Z_s, s\le t\racc$ for the running maximum of $Z,$ and recall that
$$\pb [T_x > t]\; = \;\pb[M_t \le x].$$
If $\mu$ is concentrated on $\rl^+$ then $Z$ is a.s. increasing and one is reduced to the random walk case because 
$$\pb[Z_{[t] +1} \le x]\; \le \;\pb [T_x > t]\; = \;\pb[Z_t \le x]\; = \;\pb[M_t \le x]\; \le \;\pb[Z_{[t]} \le x].$$
If $\mu$ is not concentrated on $\rl^+$ then an argument analogous to that of the previous paragraph shows that 
$$\pb [T_x > t]\; \asymp \; \pb [T_y > t]$$
as $t\to +\infty$ for every $x,y > 0,$ as in the discrete time framework. The classical approach to obtain more information on $\pb [T_x > t]$ relies on a particular case of the first fluctuation identity - see e.g. Theorem VI.5 in \cite{Be}, which is an analogue of Spitzer's formula in continuous time:

\vspace{2mm}

\noindent
{\bf Baxter-Donsker's formula.} For every $\lbd \ge 0$ and $q > 0$, one has
\begin{equation}
\label{BD}
\EE[e^{-\lbd M_{\tau_q}}]\; = \; \exp-\lcr\int_0^\infty \frac{e^{-qt}}{t}\, \EE [(1 -e^{-\lbd Z_t^+})] \,dt\rcr
\end{equation}
with the notation $Z_t^+ = \max (0, Z_t)$ and $\tau_q\sim$ Exp $(q)$ an independent random time.
\vspace{2mm}

An important consequence of this formula is a theorem of Rogozin - see e.g. Theorem VI.5 in \cite{Be} - which shows that if the following Spitzer's condition
\begin{equation}
\label{Bubi}
\frac{1}{t} \int_1^t \pb [Z_s \ge 0] ds \; \to\; \rho \in (0,1),\quad t\to +\infty
\end{equation}
holds, then for any $x > 0$
\begin{equation}
\label{Mme}
\pb [T_x > t]\; \sim\; c_x l(t) t^{-\rho} \quad \mbox{}
\end{equation}
with $c_x >0$ and $l(t)$ some slowly varying function not depending on $x$. Besides, one can show that (\ref{Bubi}) and (\ref{Mme}) are actually equivalent - see again Theorem VI.5 in \cite{Be}. Recall - see Chapter 7 in \cite{Do} - that (\ref{Bubi}) is also equivalent to (\ref{Spitz}) for the random walk $\lacc -Z_n, n\ge 1\racc$. The estimate (\ref{Mme}) can be refined for strictly $\a-$stable processes, which all enjoy the property that 
\begin{equation}
\label{Posi}
\pb [Z_t \ge 0] \; = \; \rho \in (0,1)\quad\mbox{for every $t > 0.$}
\end{equation}
An asymptotic analysis of the so-called Darling integral - see Theorem 3b in \cite{Bi1} - entails then 
\begin{equation}
\label{Bing}
\pb [T_x > t]\; \sim\; c x^{\a\rho} t^{-\rho}
\end{equation}
for such processes, where $c$ is an explicit constant. Notice that there are other L\'evy processes enjoying the property (\ref{Posi}), like subordinate stable processes, for which no refinement of (\ref{Mme}) seems available in the literature. In \cite{Ba}, the precise estimate
$$\pb [T_x > t]\; \sim\; c_x t^{-1/2}$$
was obtained for every centered L\'evy processes with finite variance and every $x >0.$ \\

We conclude this paragraph with spectrally negative L\'evy processes, where our survival analysis amounts to the study of a one-dimensional probability as in the discrete framework, and does not really require (\ref{BD}). Indeed, the passage-time process $\lacc T_x, x\ge 0\racc$ is then a subordinator with infinite lifetime if $Z$ does not drift to $-\infty$, whose Laplace exponent $\Phi(\lbd) = -\log \EE[e^{-\lbd T_1}]$ is characterized by the law of $Z_1$ - see Chapter VII in \cite{Be} or Chapter 9 in \cite{Do} for these basic facts. In particular, if $\Phi(\lbd)\sim \lbd^{\rho}l(\lbd)$  as $\lbd\to 0$ for some $\rho\in (0,1)$ and $l(\lbd)$ some slowly varying function, then Theorem XIII.5.4. in \cite{F} yields
$$\pb [T_x > t]\; \sim\; \frac{x}{\Gamma(\rho)} l(t^{-1}) t^{-\rho}.$$ 
Actually, the above condition on $\Phi$ is equivalent to (\ref{Bubi}) - see Proposition VII.6 in \cite{Be} and notice that then necessarily $\rho \ge 1/2.$ Hence, the above estimate is just a consequence of Rogozin's theorem with an explicit constant $c_x$. Spectrally negative L\'evy processes also enjoy the following peculiar property, which follows easily from the Baxter-Donsker formula expressed with the exponent $\Phi$:
\begin{equation}
\label{Skoro}
M_{\tau_q}\; \elaw\; Z_{\tau_q}\; \vert\; Z_{\tau_q} > 0.
\end{equation}
In particular, one has $M_t \,\elaw\, Z_t \,\vert\, Z_t > 0$ for every $t > 0$ if $\pb [X_t > 0] = \rho\ge 1/2$ does not depend on $t$ (which is true only for the $(1/\rho)$-stable process). This latter identity which can be shown in many different ways - see Section 8 in \cite{Bi2} and the references therein, recovers the estimate (\ref{Bing}) in this particular case: one has 
$$\pb [T_x > t]\; \sim\; \frac{ c x}{\Gamma(\rho)} t^{-\rho} \quad \mbox{}$$ 
for some explicit constant $c > 0$ which is $\sqrt{2}$ for the standard Brownian motion. We finally stress that (\ref{Skoro}) can be useful for spectrally negative processes such that $\pb [X_t > 0] \to 1.$  For example, if $X$ is an $\a-$stable process with positive drift then (\ref{Skoro}) shows after a simple analysis that
$$\pb [T_x > t]\; \sim\; c x t^{-\a}$$ 
for some explicit $c > 0,$ an estimate which cannot be obtained directly neither from Rogozin's theorem nor from the behavior of $\Phi$ at zero, and which is also coherent with (\ref{Embre}) since $\pb[X_t < 0] \sim c t^{1-\a}.$

\section{Recent advances}

\subsection{Integrated random walks} 

\label{sec:irw}

An integrated random walk is the sequence of partial sums $A_n = S_1 + \cdots + S_n$, where $\lacc S_n, \, n\ge 1\racc$ is a random walk. As above we write $S_n = X_1 + \cdots + X_n$ and we denote by $\mu$ the law of the increment $X_1 = S_1.$ Let
$$T_x \; =\; \inf\lacc n\ge 1, \; A_n > x\racc$$
be the first-passage time above $x \ge 0.$ Since $\lacc S_1\le 0, \ldots, S_n\le 0\racc \subset \lacc A_1\le 0, \ldots, A_n\le 0\racc,$ the discussion made in Paragraph \ref{sec:rw} show that $\pb[T_x = +\infty]\ge \pb[T_0 = +\infty] > 0$ as soon as $S_n$ drifts to $-\infty$. If $\mu$ is concentrated on $\rl^+,$ then $A_n$ has non-negative increments and since $A_n = nX_1 + (n-1)X_2 \cdots + X_n \ge (n/2) S_{[n/2]},$ one has
$$\pb[T_x > n]\; =\; \pb [A_n < x]\; \le \; \pb [n S_{[n/2]} < 2x]\; \le\; e^{2x + n\log (\EE[e^{-[n/2] X_1}])},$$
which shows that $\pb[T_x > n]$ tends to zero superexponentially fast, unless $\mu$ is degenerate. If $\mu$ is not concentrated on $\rl^+,$ then choosing $\eps$ such that $\mu (-\infty, -\eps) >0$ and using the same argument as in Paragraph \ref{sec:rw} entail $\pb[T_x > n] \ge \pb[T_0 > n] \ge (\mu (-\infty, -\eps))^k\pb[T_x > n]$ for every $x \ge 0$ as soon as $k\eps \ge x,$ so that 
$$\pb[T_x > n]\; \asymp\; \pb[T_0 > n]$$
for every $x\ge 0.$ \\

Henceforth we will suppose that $S_n$ does not drift to $-\infty$ and that $\mu$ is not concentrated on $\rl^+$. This entails that $\pb[S_n> 0] >0$ and $\pb[S_n< 0] >0$ for every $n\ge 1$. We are interested in the rate of decay of $\pb[T_0 > n]$ to zero, and we will see that far much less is known than for random walks.\\

The case of integrated simple random walks was first considered by Sina\u {\i} \cite{Si1}, who  showed the following

\begin{THEO} [Sina\u {\i}]Suppose that $\mu \{+1\} = \mu\{-1\} = 1/2.$ Then 
\begin{equation}
\label{quarter}
\pb[T_0 > n]\; \asymp\; n^{-1/4}.
\end{equation}
\end{THEO}

The main idea lying behind this theorem is very simple. Let $\{\tau_k, k\ge 0\}$ be the a.s. finite sequence of return times of $S_n$ to zero viz. $\tau_0 = 0$ and $\tau_k = \inf\{ n> \tau_{k-1}, \; S_n =0\}.$ On the one hand, the simplicity assumption entails that $\lacc A_{\tau_n}, n\ge 1\racc$ is an integer-valued symmetric random walk, and that the a.s. identification
$$\lacc A_{\tau_1}\le 0, \ldots, A_{\tau_n}\le 0\racc\; =\; \lacc A_i\le 0, \;\forall\; i =1\ldots \tau_n\racc$$
holds. The series $\sum_{n\ge 1} n^{-1} \pb[A_{\tau_n} = 0]$ is convergent, and Sparre-Andersen's formula combined with the Tauberian theorem for monotonic sequences shows that
$$\pb[T_0 > \tau_n]\; \sim\; \frac{c}{\sqrt{n}}\cdot$$
On the other hand, the sequence $\tau_n$ grows like $n^2$ at infinity (more precisely, $n^{-2}\tau_n$ converges in law to some positive $(1/2)-$stable law), so that after some residual analysis on the bivariate random walk $\lacc A_{\tau_n}, \tau_n\racc,$ one obtains the desired result. \\

Since then several authors have tried to extend the validity of the estimate (\ref{quarter}) to more general random walks, and the following conjecture is made in \cite{CD} and \cite{Vy0}:

\begin{Conj} Suppose that $\mu$ has finite second moment and zero mean. Then 
$$\pb[T_0 > n]\; \asymp\; n^{-1/4}.$$
\end{Conj}

A general result in this direction was obtained recently by Dembo and Gao \cite{DG}, who proved the following 

\begin{THEO}[Dembo-Gao] Suppose that $\mu$ has finite second moment and zero mean. Then there exists an explicit constant $K$ such that
$$\pb[T_0 > n]\; \le\; K n^{-1/4}, \quad n\ge 1.$$
If in addition $\exists \,\beta > 0$ such that $\mu(-\infty, -r) > e^{-\beta r}$ for some $r > 0$ and $\mu(-\infty, -t) < e^{-\beta t}$ as $t\to +\infty,$ then
$$\pb[T_0 > n]\; \asymp\; n^{-1/4}.$$ 
\end{THEO}

The method used \cite{DG} is completely different from Sina\u {\i}'s and relies on a decomposition of the integrated walk at its supremum, which is somehow reminiscent of Sparre Andersen's argument, and too involved to be discussed here in detail. This method  also allows for an elementary proof of the formula (\ref{Exact}) in the symmetric and absolutely continuous case - see Proposition 1.4 therein. Notice that the result in \cite{DG} is stated in a different manner than ours - see Theorem 1.2 and Remark 1.3 therein - and that the assumption for the lower bound is slightly less stringent. But as in our formulation, this assumption essentially means that $\EE [e^{-\beta X_1}] < +\infty$ for some $\beta>0$ and that the left tail $\mu(-\infty, -x)$ has some regularity in the neighbourhood of zero - see Remark 1.1 in \cite{DG}. Shortly before this result, it was shown in \cite{AD} in using strong approximation that
$$n^{-1/4}(\log n)^{-4}\; \le\; \pb[T_0 > n]\; \le\; n^{-1/4}(\log n)^4$$ 
at infinity when $\EE [e^{\beta |X_1|}] < +\infty$ for some $\beta>0$. \\

Sina\u {\i}'s method of decomposing the integrated path along the excursions away from zero of the underlying random walk can be extended under some extra assumptions on the law  of $-X_1$ conditioned to be positive, which we denote by $\mu^-$. This idea had been used in \cite{Vy1} to show the same result as \cite{DG} in several situations like double-sided exponential, double-sided geometric, left-continuous, or lazy simple random walks, which all satisfy Dembo and Gao's assumption for the lower bound. Recently in \cite{Vy2}, this technique is combined with local limit theorems to obtain the following more precise result. 

\begin{THEO}[Vysotsky] Suppose that $\mu_-$ is exponential or that $\mu$ is left-continuous. If $\mu$ has finite variance and zero mean then there exists $c = c(\mu) >0$ such that
$$\pb[T_0 > n]\; \sim \; c \,n^{-1/4}.$$ 
\end{THEO}

Up to now, Conjecture 1 is still unsolved in general, and the main challenge is to get rid of the extra assumptions on $\mu_-.$ As can be seen from the central limit theorem, random walks with zero mean and finite variance are such that $\pb[S_n > 0]\to 1/2.$ In view of the discussion made in Paragraph \ref{sec:rw}, it is hence natural to raise the more general

\begin{Conj} Suppose that $\pb[S_n > 0]\to 1/2.$ Then 
$$\pb[T_0 > n]\; = \; n^{-1/4 +o(1)}.$$
\end{Conj}

In this formulation, no assumption is made on the moments of $\mu$ and this enhances sharply the difficulty of the problem, which probably requires more combinatorial tools than the one used in the above references. \\

We now turn to some situations where the persistence exponent of integrated random walks is not 1/4. We denote by $\cD(\a)$ the set of probability measures attracted to some strictly $\a-$stable law with $\a\neq 2$ and we refer e.g. to Chapter XVII.5 in \cite{F} for more on the subject. Recall that if $\mu\in\cD(\a),$ then $\EE[\vert X_1\vert^s ] < \infty$ for every $s\in [0, \a)$ and that $\EE[X_1] = 0$ if $\a > 1.$ The following is a consequence of the main result in \cite{DG}:

\begin{THEO}[Dembo-Gao] Suppose that $\mu\in\cD(\a)$ for some $\a\in (1,2).$ Then there exists an explicit constant $K$ such that
$$\pb[T_0 > n]\; \le\; K n^{-(1-1/\a)/2}, \quad n\ge 1.$$
If in addition $\mu$ is attracted to a spectrally positive $\a-$stable law, then
$$\pb[T_0 > n]\; \asymp\; n^{-(1-1/\a)/2}.$$ 
\end{THEO}
The result in \cite{DG} is formulated in a different manner and is actually more general. In the case when $\mu\in\cD(\a)$ however, the additional assumption made therein for the lower bound is equivalent to the spectral positivity of the attracting law - see e.g. Theorem XVII.5.1 in \cite{F}. With the help of Sina\u {\i}'s method, the estimate was recently refined in \cite{Vy2}:

\begin{THEO}[Vysotsky] Suppose that $\mu_-$ is exponential or that $\mu$ is left-continuous. If $\mu$ is normally attracted to some spectrally positive $\a-$stable law with $\a\in (1,2)$ then there exists $c = c(\mu) >0$ such that
$$\pb[T_0 > n]\; \sim \; c \,n^{-(1-1/\a)/2}.$$ 
\end{THEO}

The result in \cite{Vy2} is stated differently but we see from \cite{BGT} p. 382 that it is actually equivalent to our formulation. In view of our final discussion made in Paragraph \ref{sec:rw}, it is very surprising that the case where $\mu$ is right-continuous and normally attracted to some spectrally negative $\a-$stable law seems to be more difficult to handle than the dual situation where $\mu$ is left-continuous and normally attracted to some spectrally positive $\a-$stable law.\\

It is natural to ask what the persistence exponent should be when the limit stable law has negative jumps. Let us hence denote by $\cD(\a, \rho)$ the set of probability measures attracted to a strictly $\a-$stable law with positivity parameter $\rho.$ Notice - see e.g. \cite{Z} for details and recall that we excluded the one-sided case - that $\rho\in (0,1)$ for $\a\in (0,1]$, that $\rho\in [1-1/\a,1/\a]$ for $\a\in (1,2),$ and that spectrally positive $\a-$stable laws with $\a\in (1,2)$ are such that $\rho = 1-1/\a.$ Besides, one has $\pb[S_n > 0] \to \rho$ whenever $\mu\in\cD(\a, \rho)$ - see e.g. Theorem XVII.5.1 in \cite{F}. In view of the discussion made in Paragraph \ref{sec:rw}, it is natural to raise the general

\begin{Conj} Suppose that $\pb[S_n > 0] \to \rho\in(0,1).$  Then 
$$\pb[T_0 > n]\; =\; n^{-\rho/2 +o(1)}.$$
\end{Conj}
A weak version of the conjecture is to get the asymptotic under the assumption that $\mu\in\cD(\a, \rho)$, just like Conjecture 1 which deals with the case $\a = 2$ and $\rho = 1/2.$ Under the stronger assumption that $\mu$ is normally attracted, one may also wonder if a more precise behavior could not be obtained, as in \cite{Vy2}.

\subsection{Integrated L\'evy processes} 

\label{sec:ilp}

In this section we consider the process
$$A_t \; =\; \int_0^t Z_s \, \d s, \quad t\ge 0,$$
where $\lacc Z_t, \; t\ge 0\racc$ is a real L\'evy process starting from zero. We set 
$$T_x \; = \; \inf\lacc t > 0, \; A_t > x\racc\; = \; \inf\lacc t > 0, \; A_t = x\racc$$
for its first passage time above $x\ge 0.$ Contrary to (\ref{Rough}), there is no straightforward  bound between $\pb[T_x \ge t]$ and an analogous quantity involving some iterated random walk, so that the results of the previous paragraph cannot be used directly. 

The process $(Z,A)$ is strongly Markovian and we set $\pb_{(z,a)}$ for its law starting from $(z,a),$ with the simplified notation $\pb = \pb_{(0,0)}.$ It is clear by the right-continuity of $Z$ that 
$$\pb_{(z,0)}[T_0 = 0] = 0 \;\rm{or}\; 1 \quad\mbox{according as $z >0$ or $z<0.$}$$ 
One has also $\pb[T_0 = 0] = 0$ or $1$ by the 0-1 law, but to obtain a criterion for the regularity of the upper half-plane for $(Z,A)$ is an open problem which does not seem obvious. Since integrated L\'evy process all have finite variation, one might think that this criterion will be different from the aforementioned Theorem 22 in \cite{Do}.

If $Z$ drifts to $+\infty$, then it is clear that $T_x < +\infty$ a.s. for every $x > 0.$ On the other hand, if $Z$ drifts to $+\infty$ then its last passage time above zero can be made arbitrarily small so that one will have $\pb[T_x = +\infty] > 0$ for every $x > 0.$ When $Z$ oscillates then probably one has $\pb[T_x = +\infty] = 0$ for every $x > 0,$ but there is not a direct answer to this question. In general there is no result of basic fluctuation theory available for integrated L\'evy processes. \\

In this paragraph we will consider two examples where the persistence exponent can be computed. The first one is the integrated Brownian motion and originates from Kolmogoroff \cite{Ko}, in relation with the two-dimensional generator
$$\frac{1}{2}\, \frac{\partial^2}{\partial x^2}\; +\; x \, \frac{\partial}{\partial y}$$
and the associated Fokker-Planck equation. Notice that \cite{Ko} actually deals with the more general $n$-times integrated Brownian motion. The process $(B,A)$ is a Gaussian Markov process whose transition density can be computed via the covariance matrix: under $\pb = \pb_{(0,0)}$ one has
$$p_t(b,a) \; =\; \frac{\sqrt{3}}{\pi t^2}\exp[- 2b^2/t + 6ab/t^2 -6 a^2/t^3],$$
and the expression under $\pb_{(x,y)}$ follows by translation. By the (3/2)-self-similarity of $A$, one has
$$T_x \;\elaw\; x^{2/3} T_1$$
under $\pb,$ so that our persistence problem amounts to find the asymptotic of $\pb[T_1 > t],$ a question which dates back to Uhlenbeck and Wang in 1945. Notice that the above identity also yields $\pb[T_0 = 0] = 1.$ Among other formul\ae, the following was obtained by McKean in an analytical way - see (3.1) in \cite{McK}:
$$\pb_{(-1,0)}\lcr T_0 \in dt, B_{T_0} \in dx\rcr\; =\; \frac{3x}{\pi \sqrt{2\pi} t^2}e^{- (2/t)(1-x+x^2)}\lpa\int_0^{4x/t} \!\!\!\!e^{-3y/2}y^{-1/2} dy\rpa\Un_{\{x\ge 0\}} dt dx.$$
This formula is the key argument to the following result which is proved separately in \cite{G} and \cite{IW}:

\begin{THEO}[Goldman, Isozaki-Watanabe] With the above notations there exists $c > 0$ such that
\begin{equation}
\label{GIW}
\pb[T_1 > t]\; \sim\; c\, t^{-1/4}.
\end{equation}
\end{THEO}
\noindent
Notice that by self-similarity this result has the more general formulation 
$$\pb [A_s\le x, \;\forall\, s \in [0,t]] \; = \; \pb[T_x > t]\; \sim\; c\, x^{1/6}t^{-1/4} \quad \mbox{as $t/x^{2/3}\to +\infty.$}$$
In particular, one has
$$\pb [A_s\le x, \;\forall\, s \in [0,1]] \;\sim\; c\, x^{1/6} \quad \mbox{as $x\to 0,$}$$
which is a lower tail probability statement as mentioned in the introduction. By strong approximation and as a consequence of his result on integrated simple random walks, Sina\u {\i} had obtained in \cite{Si1} the rougher estimate
$$\pb[T_1 > t]\; \asymp\; t^{-1/4}.$$
The argument in \cite{G} is analytical and relies on integral equations. It also yields a complicated explicit expression for the law of $T_1:$ 
\begin{eqnarray*}
\pb\lcr T_1 \in dt \rcr & = &\pb_{(0,-1)}\lcr T_0 \in dt \rcr\; =\; \lpa\frac{3\sqrt{3}}{2\sqrt{2\pi}t^{5/2}}e^{-3/2t^3}\right. \\
 &+ &\!\!\!\!\left.\frac{2\sqrt{3}}{\pi}\!\!\int_0^\infty\!\!\!\!\int_0^\infty\!\!\!\!\int_0^t \!\!\pb_{(-z,0)}\lcr t -T_0 \in ds, B_{T_0} \in dx\rcr e^{-6/s^3-2z^2/s^2}\sinh(6z/s^2)  zdzds\rpa \! dt,
\end{eqnarray*}
where the quantity $\pb_{(-z,0)}\lcr t -T_0 \in ds, B_{T_0} \in dx\rcr$ can be expressed via McKean's above formula and self-similarity. Notice that both McKean and Goldman's formul\ae \, have been generalized by Lachal \cite{La}, who obtained an explicit formula for $\pb_{(b,a)}\lcr T_0 \in dt, B_{T_0} \in dx\rcr$ for any $(b,a).$ The asymptotic analysis which is carried out in Proposition 2 of \cite{G} gives the right speed of convergence for $\pb[T\in \d t]$, but the resulting value of the constant $c$ in (\ref{GIW}) is erroneous because of some inaccuracies in the change of variable. The right value is
$$c\; =\; \frac{3^{4/3}\Gamma(2/3)}{\pi 2^{13/12} \Gamma(3/4)}\;\sim\; 0.718$$
and follows also from the simpler probabilistic method of \cite{IW}, which relies on the Markov property and a Tauberian argument. This method yields the more general estimate 
$$\pb_{(b,0)}[T_x > t]\; \sim\; c_{b,x}t^{-1/4}$$
with some explicit $c_{b,x} >0$ - see (1.12) in \cite{IW}. It also allows to handle first-passage time asymptotics for fluctuating homogeneous additive functionals of Brownian motion, with a persistence exponent depending smoothly on the skewness of the functional - see Corollary 1 in \cite{IK}. We finally notice that the estimate (\ref{GIW}) was also obtained (with an erroneous constant $c$) in \cite{Bu} after solving some Krein-Kramers differential equation in the context of semiflexible polymers in the half-plane. This latter method was generalised in \cite{BG} to give another computation of the persistence exponent for fluctuating homogeneous additive functionals of Brownian motion, in the context of survival of a diffusing particle in a transverse shear flow. \\

We now turn to integrated strictly $\a$-stable L\'evy processes, which form the natural generalisation of integrated Brownian motion. The bivariate process $(Z,A)$ is then a stable Markov process, where the stability property has the general meaning which is given in the monograph \cite{ST}. In particular, one can check from the L\'evy-Khintchine formula - see e.g. Proposition 3.4.1 in \cite{ST} - that 
$$A_1 \;\elaw\; (1+\a)^{-1} Z_1$$
for every $t \ge 0,$ so that $A_1$ is a $\a-$stable variable with the same positivity parameter as $Z_1.$ On the other hand there does not seem to exist an explicit formula for the density of the bivariate random variable $(Z_1, A_1)$ except in the Gaussian case $\a = 2.$ The univariate process $A$ is $(1+1/\a)-$self-similar, so that with the above notations one has
$$T_x \;\elaw\; x^{\a/(\a +1)} T_1$$
under $\pb$ and we need only to study the asymptotic behaviour of $\pb[T_1 > t].$ Notice that again this identity yields $\pb[T_0 = 0] = 1.$ As a consequence of the main result of \cite{TS1}, the persistence exponent of $A$ can be computed within a specific sub-class:

\begin{THEO}[Simon] Let $\lacc Z_t, \; t \ge 0\racc$ be a strictly $\a$-stable L\'evy process with $\a\in (1,2).$ With the above notations, there exists a positive constant $K$ such that
$$\pb[T_1 > t]\; \le\; K\, t^{-(1-1/\a)/2}, \quad t > 0.$$
If in addition $Z$ is spectrally positive, then 
$$\pb[T_1 > t]\; =\; t^{-(1-1/\a)/2 +o(1)}.$$ 
\end{THEO}

The main result in \cite{TS1} deals with more general homogeneous functionals of stable L\'evy processes, extending the results of \cite{Is}. It also provides some explicit lower bound with a logarithmic term, which entails the following criterion for the finiteness of fractional moments of $T_1$ in the spectrally positive case:
$$\EE [T_1^s] < \infty\; \Leftrightarrow\; -(\a + 1) < s < (\a-1)/2\a.$$
The method of \cite{TS1} is an adaptation of Sina\u {\i}'s in continuous time, relying on the bivariate L\'evy process $\lacc (\tau_t, A_{\tau_t}), \; t\ge 0\racc$ with $\lacc \tau_t, \, t\ge 0\racc$ the inverse local time at zero (which exists because $\alpha > 1$). The upper bound relies on the Wiener-Hopf factorization method as in \cite{Is}, and the fact that $A_{\tau_1}$ is a symmetric $(\a -1)/(\a+1)$-stable variable whatever $A_1$'s positivity parameter is. For the lower bound, the crucial fact is that a.s.
$$\lacc A_{\tau_s} \le 1,\; \forall s\in[0,t]\racc\; =\; \lacc A_s \le 1,\; \forall s\in[0,\tau_t]\racc,$$
which allows to study the probability of the right event with the help of (\ref{Bing}). This latter identity is true only in the spectrally positive case. \\

It is a natural question to find the persistence exponent for all integrated stable L\'evy processes. If $Z$ is an $\a$-stable subordinator, then $A$ is an increasing process and one has
$$\pb[T_1 > t]\; =\; \pb[A_1 \le t^{-(1+1/\a)}]\; \sim\; c_\a t^{-(1+\a)/2\a(1-\a)}e^{-(1-\a)\a^{\a/(1-\a)}t^{-(1+\a)/(1-\a)}}$$   
where the estimate follows e.g. from (2.4.30) in \cite{IL}. This superexponential speed is coherent with the previous discussion for integrated one-sided random walks. If $-Z$ is an $\a$-stable subordinator, then $T_x = +\infty$ a.s. for every $x\ge 0.$ The general problem remains open, and in view of Bingham's estimate (\ref{Bing}) it is natural to raise the precise

\begin{Conj} Let $Z$ be a strictly $\a-$stable L\'evy process such that $\pb[Z_1 > 0] =\rho\in(0,1).$ Then there exists $c >0$ such that
$$\pb[T_1 > t]\; \sim\; c\, t^{-\rho/2}.$$
\end{Conj}
We see that the main result of \cite{TS1} shows the conjecture in the spectrally positive case viz. $\a > 1$ and $\rho = 1- 1/\a$, except the existence of the constant for which the refined methods of \cite{Vy2} could perhaps be adapted. In the general case, a simple scenario with an initial downward intention is described in \cite{DST} p. 3 to show that the speed must be larger than $t^{-\rho}.$ Unfortunately, the attempts made therein to obtain a better bound were unjustified and in \cite{TS2} a method relying on decorrelation inequalities, inspired by \cite{CD}, is carried out to show that the upper bound in the main result of \cite{TS1} is not the optimal one when $Z$ has negative jumps and self-similarity index close to 1:

\begin{THEO} [Simon] For every $\a_0, c >0$ there exists $k > 0$ such that for every $\a\in[\a_0,2]$ and every strictly $\a-$stable process $Z$ with positivity parameter $\rho\in[c\wedge (1-1/\a), 1/\a\vee 1]$ one has
$$\liminf_{t\to \infty} t^k\pb[T_1 > t]\; =\; 0.$$
\end{THEO}
At present, it does not seem that the different arguments developed in \cite{DST, TS1, TS2}  can provide any clue to get the conjectured persistence exponent $\rho/2.$ It is quite surprising that when $Z$ is spectrally positive the persistence exponent is easier to get for $A$ than for $-A$. Indeed, the simple identity (\ref{Skoro}) entails readily that the persistence exponent of $-Z$ is $1/\a,$ whereas it is not a straightforward task to show that the persistence exponent of $Z$ is $1-1/\a.$ \\

It seems plausible that a version without constant of Conjecture 4 could be proved in refining the results of \cite{DG} and using some discrete approximation. In view of  Conjecture 3 and (\ref{Mme}), this leads us to the following general question on integrated L\'evy processes:

\begin{Conj} Let $Z$ be a L\'evy process such that $\pb[Z_t > 0] \to \rho\in(0,1)$ as $t\to +\infty.$ Then for every $x >0$
$$\pb[T_x > t]\; =\; t^{-\rho/2 +o(1)}.$$
\end{Conj}

\subsection{Fractionally integrated L\'evy processes} 

\label{sec:filp}

In this section we consider processes of the type
\begin{equation}
\label{RLP}
\Ab_t \; =\; \frac{1}{\Gamma(\b+1)}\int_0^t (t-s)^\b \, \d Z_s, \quad t\ge 0,
\end{equation}
where $\lacc Z_t, \; t\ge 0\racc$ is a real L\'evy process starting from zero and $\b > -1.$ The above convolution product makes sense for every $\b \ge 0,$ and can also be defined for some negative $\b$ depending on the law of $Z.$ If $Z$ is strictly $\a-$stable for instance, then $\Ab$ is well-defined for every $\b > -(1\wedge1/\a)$ and is then a stable process in the sense of \cite{ST}, with continuous paths iff $\a = 2$ or $\b > 0$ and with a.s. locally unbounded paths iff $\a < 2$ and $\b < 0$ - see Chapter 10 in \cite{ST}. When $Z$ is strictly $\a-$stable, it is customary to write $\b = H - 1/\a$ with $H >0$ the so-called Hurst parameter, and $\Ab$ is then an $H-$self-similar process.  One can view $\b$-fractionally integrated L\'evy processes as the natural generalization of $n$ times integrated L\'evy processes, since an integration by parts shows that the latter form the subclass $\b = n.$ In this particular case the $(n+1)$-dimensional process $(Z, A^1, \ldots, A^n)$ is Feller, but there is no multidimensional Markov property when $\b$ is not an integer because then the fractional integration takes the whole memory of the driving process into account. In the literature fractionally integrated L\'evy processes are often called Riemann-Liouville processes, a denomination which is originally due to L\'evy. 

In the Brownian case $Z = B,$ the process $\Ab$ is closely connected to the fractional Brownian motion $\lacc B^H_t, \; t\ge 0\racc$, which we recall to be the centered Gaussian process with covariance function
$$\EE \lcr B^H_t B^H_s\rcr\; =\; \frac{1}{2} \lpa t^{2H} + s^{2H} - |t-s|^{2H}\rpa, \qquad t,s\ge 0.$$
Fractional Brownian motion can be written as the independent sum 
\begin{equation}
B^H_t\; =\; c_H\lpa A^{H-1/2}_t\; +\; \int_0^{\infty}((t+s)^{H-1/2} - s^{H-1/2}) \d {\tilde B}_s\rpa, \quad t\ge 0,
\end{equation}
with $c_H = (H2^{2H}\pi/\Gamma(H+1/2)\Gamma(1-H))^{1/2}$ the normalization constant, which shows that its paths are continuous a.s. The process $B^H$ gives insight on the long-range increments of $A^{H-1/2}$ since it can be proved in analysing the covariance function that
\begin{equation}
\label{Incr}
\{ A^{H-1/2}_{t+u} - A^{H-1/2}_u, \; t\ge 0\}\; \claw\; \lacc c_H^{-1}B^{H}_{t}, \; t\ge 0\racc, \quad u\to +\infty
\end{equation}
(with, of course, an equality in law for every $u$ when $H =1/2$). Recall that $B^H$ is defined for $H\in (0,1]$ only and that $B^1$ is simply the linear function $t\mapsto t N$ with $N$ a standard normal variable. Fractional Brownian motion can be shown \cite{ST} to be the unique $H-$self-similar Gaussian process with stationary increments, whence its greater importance in modeling than Riemann-Liouville processes. Setting 
$$T^H_x\; =\; \inf\lacc t > 0, \; B^H_t > x\racc \; = \;\inf\lacc t > 0, \; B^H_t = x\racc$$ for every $x\ge 0,$ one has $T^H_x\, \elaw\, x^{1/H} T_1^H$ by self-similarity, which shows that $T^H_0 = 0$ a.s. and that the survival analysis of $B^H$ is reduced to the behavior of $\pb[T^H_1 > t]$ only. Among other results, the following was obtained in \cite{Mo1}:

\begin{THEO}[Molchan] For every $H\in(0,1]$ one has
$$\pb[T^H_1 > t]\; =\; t^{H-1+o(1)}.$$
\end{THEO}
\noindent
The main argument of \cite{Mo1}, partly inspired by Brownian fluctuation theory, is to quantify the correlation between $T^H_1,$ the last zero of $B^H$ on $[0,1],$ the positive sojourn time of $B^H$ on $[0,1],$ and the inverse exponential functional
$$t \; \mapsto\; J^H_t\; =\;\lpa \int_0^t e^{B^H_s} \d s \rpa^{-1},$$
whose asymptotic behavior in expectation can be precisely analysed. The general link between $J^H_t$ and $T^H_1$ is explained by the heuristical fact that if $T^H_1 > t$ then $B^H$ has drifted towards $-\infty$ rather rapidly, so that $J^H_t$ is big. Conversely if $T^H_1 < t$ then $B^H$ has been close zero for a positive fraction of time, so that $J^H_t$ is small. However, the analysis of $\{\EE[J^H_t], \; t\ge 0\}$ which is performed in \cite{Mo1} is very specific to the stationary increments of fractional Brownian motion, and does not seem to be suitable for Riemann-Liouville processes. \\

It is clear that $\pb[T^1_1 > t]\to 1/2$ by the above remark on the case $H= 1$, and it is natural to raise the 
\begin{Conj} For every $H\in (0,1)$ one has 
$$\pb[T^H_1 > t]\; \asymp\; t^{H-1}.$$
\end{Conj}
In view of the previously stated precise results, one may also ask if $\pb[T^H_1 > t] \sim c\, t^{H-1},$ but from empirical studies which were communicated to us by A.~A.~Novikov, it seems that such an exact behaviour does not hold in general. According to these studies, at least for $H= 3/4$ one should have
 $$\liminf_{t\to +\infty} t^{1-H}\pb[T^H_1 > t]\; < \; \limsup_{t\to +\infty} t^{1-H}\pb[T^H_1 > t].$$

The above conjecture was recently addressed in \cite{A}, where the following partial result is obtained by a refinement of Molchan's method:

\begin{THEO}[Aurzada] There exists $c>0$ such that
$$(\log t)^{-c}t^{H-1}\; \le\;\pb[T^H_1 > t]\; \le\;(\log t)^ct^{H-1}, \quad t\to +\infty.$$
\end{THEO}

We now turn to first passage asymptotics for a class of Gaussian stationary processes (GSP) which are related to the persistence of fractionally integrated Brownian motion. In \cite{LS4} the Lamperti transformation
\begin{equation}
\label{Lamp}
L^H_t\; =\; e^{-tH}B^H_{e^t}, \quad t\in\rl,
\end{equation}
a centered GSP, is studied in connection with the persistence exponent of $B^H$ and it is shown that
$$\log \pb\lcr L^H_s \le 0,\; \forall\, s\in[0,t]\rcr\; =\; \log \pb\lcr B^H_s \le 0,\; \forall\, s\in[1,e^t]\rcr\;\sim\; t (H-1).$$  
Since Molchan's theorem means that $\log \pb\lcr B^H_s \le 1,\; \forall\, s\in[0,e^t]\rcr\,\sim\, t (H-1),$ one simply needs to switch the 0 and the 1 in the latter probability to obtain this result, which is justified in \cite{LS4} through a refined use of Slepian's lemma. The Lamperti transformation is also a fruitful method in the reverse direction, and this was observed in \cite{AD} to investigate the persistence exponents of Riemann-Liouville processes. Set $Z =B$ in (\ref{RLP}) and introduce the notation
$$T^\b_x \; =\; \inf\lacc t > 0, \; A^\b_t > x\racc\; =\;  \inf\lacc t > 0, \; A^\b_t = x\racc$$
for every $\b, x \ge 0.$ Notice that $T^\b_x \,\elaw\, x^{2/(1+2\b)}T^\b_1$ by self-similarity, which shows that $T^\b_0 = 0$ a.s. and that the survival analysis of $B^H$ is reduced to the behavior of $\pb[T^\b_1 > t]$ only. The process 
\begin{equation}
\label{Lampada}
Y^H_t\; =\; e^{-tH}A^{H-1/2}_{e^t}, \quad t\in\rl,
\end{equation}
a centered GSP for every $H > 0$ which reduces to the stationary Ornstein-Uhlenbeck process when $H = 1/2,$ plays a key r\^ole in the following 
\begin{THEO}[Aurzada-Dereich] There exists a non-increasing function $\b \mapsto\theta(\b)$ such that
$$\pb[T^\b_1 > t]\; =\; t^{-\theta(\b) +o(1)}$$
for every $\b\ge 0.$ Besides one has $\theta (\infty) = - \lim_{t\to\infty}t^{-1}\log \pb\lcr Y_s < 0,\; \forall\, s\in[0,t] \rcr,$ where $Y$ is the centered {\em GSP} with correlation $\EE [Y_0 Y_t] = 1/\cosh (t/2).$
\end{THEO}
The important process $Y$ is mentioned in \cite{M} in the context of diffusion equation with white noise initial conditions, and in \cite{DPSZ} in relation with the positivity of random polynomials with large even degree. More details will be given in the next section.  In \cite{LS4}, the process $Y$ is also viewed as the Lamperti transformation of the curious smooth (1/2)-self-similar Gaussian process
$$X_t\; =\; \sqrt{2}t^2\int_0^\infty B_u e^{-ut} du, \quad t\ge 0,$$
which shares the same time inversion property $\{X_t, \, t >0\}\elaw \{t^{-1}X_{t^{-1}},\, t >0\}$ as Brownian motion. In \cite{AD} only the upper bound 
\begin{equation}
\label{UB}
- \lim_{t\to\infty}t^{-1}\log \pb\lcr Y_s < 0,\; \forall\, s\in[0,t] \rcr\; \le\; \theta (\infty)
\end{equation}
is proved via Slepian's lemma, but one can easily show that $Y^H\,\claw\, Y$ as $H\to\infty$ in analysing the covariance function, so that the inequality in (\ref{UB}) is actually an equality, the limit on the left-hand side being a supremum. \\

The above result has several interesting consequences. First, it shows that the persistence exponent is a non-increasing function of the order of integration for Brownian Riemann-Liouville processes. The fact that smoother processes have more probability to survive is believed to be a kind of universal feature.  Second, it entails that $\theta (\b)\ge \theta(1) = 1/4$ for every $\b \le 1,$ so that by Molchan's result the persistence exponents of $B^H$ and $A^{H-1/2}$ do not coincide whenever $H>3/4.$ We actually believe that $\theta(H-1/2) > 1-H$ for every $H\in (1/2,1)$ and some reasons for that will be given soon afterwards. Last, it entails that $\theta (\infty)\le \theta(1) = 0.25,$ which improves the bound $\theta (\infty) < 0.325$ obtained in \cite{LS4} via another Slepian's inequality. 
In \cite{M} the numerical value $\theta(\infty) \sim 0.1875$ is suggested, whereas in \cite{DPSZ} the value $\theta(\infty) = 0.19 \pm 0.01$ is obtained by simulations. It is a tantalizing question to compute the function $\theta(\b)$ for every positive $\b\not\in\{ 0,1\},$ as well as its limit $\theta(\infty).$ The lower bound $\theta (\infty) > 0.125$ is obtained in \cite{LS2} with the help of a certain Gaussian comparison inequality, and in \cite{Mo3} this lower bound is improved into $\theta (\infty) > 1/4\sqrt{3} > 0.144,$ in comparing $Y$ with a linear time-change of the so-called Wong process \cite{W}, which is the GSP associated with integrated Brownian motion. \\

Let us give some more details on the correlation function $C_H(t) = \EE [Y^H_0 Y^H_t]$ of the Lamperti transform introduced in (\ref{Lampada}). It is given by
$$C_H(t)\; = \;2He^{-Ht}\int_0^1 (u(u + e^t -1))^{H-1/2} \d u$$
for every $t\ge 0,$ and a simple analysis shows that for every $H\in (0,1]$ and $t\ge 0$ one has
$$C_H(t)\; \le \; \cosh (Ht)\; -\; 2^{2H-1}(\sinh (t/2))^{2H}\; =\; \EE [L^H_0 L^H_t],$$
where $L^H$ is the Lamperti transform defined in (\ref{Lamp}). By Slepian's lemma and Molchan's theorem, this entails
\begin{equation}
\label{Ineq}
\theta(H-1/2)\;\ge\; 1-H
\end{equation}
for every $H \in [1/2, 1],$ with an equality if $H = 1/2.$ From the above we know that (\ref{Ineq}) is strict if $H > 3/4.$ It is very likely that $\theta(H-1/2)$ exists for every $H\in (0,1/2)$ and that (\ref{Ineq}) is strict for every $H\neq 1/2.$ In view of the above discussion, it is natural to raise the

\begin{Conj} The function $\beta\mapsto\theta(\b)$ is convex decreasing.
\end{Conj}

It is also interesting to look at the first-order expansion at zero
$$C_H(t)\; =\; 1\;-\; \frac{\Gamma(1-H)\Gamma(1/2 +H)}{\sqrt{\pi}}\lpa \frac{t}{2}\rpa^{2H} +\; o(t^{2H}),$$ 
which should be compared with 
$$\EE [L^H_0 L^H_t]\; = \; 1\;-\; \lpa \frac{t}{2}\rpa^{2H} +\; o(t^{2H}).$$ 
One might actually wonder if, for every $H\in (0,1/2)\cup(1/2,1),$ the sole inequality 
$$ \frac{\Gamma(1-H)\Gamma(1/2 +H)}{\sqrt{\pi}}\; > \; 1$$ 
could not be enough to prove that $\theta (H-1/2) > 1-H.$ Indeed, even though the asymptotic analysis of large excursions of continuous GSP's is known to be a hard problem in general, which should involve the whole correlation structure, there are also reasons to believe that its behavior at zero plays a prominent r\^ole when comparing self-similar processes having the same sample path regularity. For $H=1$ one has 
$$C_1(t)\; =\; 1\;+\; \frac{t^2}{4}\log t \; +\; o(t^2)$$ 
and from this logarithmic behaviour, which is related to the local growth of Brownian motion, one might view half-integrated Brownian motion as a kind of boundary object, analogous to the Cauchy process among stable L\'evy processes, and where the survival analysis should be different than in the other cases. This boundary phenomenon has actually been observed by physicists for a while, and we will give more details in the next section. For $H > 1$ one has
$$C_H(t)\; =\; 1\;-\; \frac{H t^2}{8(H-1)}\; +\; o(t^2),$$ 
which entails by the It\^o-Rice formula \cite{I} that the zero-crossings of the process $Y^H$ are isolated points, whose number $N^H_t$ on an interval of length $t$ has expectation
$$\EE[N^H_t]\; =\; \frac{t}{2\pi}\sqrt{\frac{H}{(H-1)}}\cdot$$ 
The interesting fact that $H\mapsto \EE[N^H_t]$ is decreasing might well be related to the non-increasingness of $H\mapsto \theta(H-1/2) = -\lim_{t\to +\infty} t^{-1} \log\pb[N^H_t =0].$ Although no mathematical result seems available, all empirical results in the physical literature \cite{M, MS} show namely that for smooth centered GSP's with positive correlation function $F$ and zero-crossing number $N_t,$ the greater $-F''(0),$ the greater $\theta = -\lim_{t\to +\infty} t^{-1} \log\pb[N_t =0].$ \\

Let us shortly notice that integrated Brownian motion is somehow critical as far as the small values of the length of the excursions of the associated Lamperti process are concerned. As a consequence of the Longuet-Higgins formula \cite{W} and a second-order expansion of $C_H(t)$ at zero, one has indeed
$$\lim_{t\to 0}t^{-1}\pb[N^H_t \neq 0]\; = \; \lacc\begin{array}{ll} +\infty & \mbox{if $H < 3/2$}\\
37/96 & \mbox{if $H = 3/2$}\\
0 & \mbox{if $H > 3/2.$}
\end{array}\right.$$
However, this has probably less to deal with the fact that $\theta(H-1/2)$ is computable only for $H = 3/2$, since analysing small and large values of the length of the excursions of continuous GSP's are very different problems. \\

We now go back to general fractionally integrated L\'evy processes and state the following result of \cite{AD}, which is obtained by strong approximation and drift transformations of Gaussian processes:

\begin{THEO}[Aurzada-Dereich] For every $\b\ge 0$, the persistence exponent $\theta(\b)$ is the same among all $\b-$fractionally integrated L\'evy processes such that $t\mapsto\EE[e^{t Z_1}]$ is defined in an open neighbourhood of zero.
\end{THEO}

The main result of \cite{AD}, to which we refer for details, is more general and handles fractionally integrated random walks as well as more general Volterra kernels. The finiteness of exponential moments play an important r\^ole therein due to the use of strong approximation, but in view of the aforementioned results and conjectures on integrated random walks and L\'evy processes, it is natural to raise the 

\begin{Conj} For every $\b\ge 0$, the exponent $\theta(\b)$ is the same among all $\b-$fractionally integrated L\'evy processes with finite variance.
\end{Conj}

We conclude this paragraph with the integrated fractional Brownian motion
$$I^H_t \; =\; \int_0^t B^H_s\, ds, \quad t\ge 0.$$
This is a centered Gaussian process with positive correlation function, and with the help of its Lamperti transformation it is easily shown that for every $H\in (0,1)$ there exists $\rho (H) > 0$ such that
$$\pb\lcr I^H_s \le 1, \; \forall\, s\in[0,t]\rcr \; =\; t^{-\rho(H) +o (1)}.$$
Numerical simulations \cite{KM} suggest the following: 

\begin{Conj}[Khokhlov-Molchan] One has $\rho (H) = H(1-H).$
\end{Conj}

This expected value, which is symmetric with respect to $H = 1/2,$ is very surprising because it is known that fractional Brownian motions with Hurst index smaller or greater than 1/2 are very different processes from several viewpoints. This does not match either the above heuristic discussion on smooth GSP's since the correlation function $F^H$ of the Lamperti process of $I^H$ has a first-order expansion at zero
$$F^H(t)\; = \; 1\; -\; \frac{(1-H^2)t^2}{2}\; +\; o(t^2).$$
In particular, the number $N^H_t$ of zero-crossings of this Lamperti process has an expectation $\EE [N^H_t] = t\sqrt{1-H^2}/\pi$ which decreases with $H$, and one might think that $\rho (H)$ also decreases. However integrated fractional Brownian motion may well be more the exception than the rule for this kind of questions, because of its complicated correlation structure. In \cite{KM} it is argued that the difference between $H <1/2$ and $H > 1/2$ should be observed at the logarithmic level and it is shown in \cite{Mo2}, with a detailed analysis, that $ c H (1-H)  \le \rho (H) \le 1-H$ for some $c\in (0,1).$ Other bounds such as $\rho(H)\ge \rho(1-H)$ for every $H\le 1/2$ were also recently presented in \cite{Mo3} (with numerical explanations), neither proving nor disproving the above conjecture.

\subsection{Other processes} Integrated random walks which were considered previously can be written as the weighted sum
$$A_n\; =\; \sum_{i=1}^n (n+1-i)X_i,$$
where $\lacc X_i, \, i\ge 1\racc$ is an i.i.d. sequence. The above weights depend on $i$ and $n$ so that the increments of $A_n$ are neither stationary nor independent. It is natural to consider persistence problem for other weighted sums like
$$\Sigma_n\; =\; \sum_{i=1}^n \sigma_i X_i,$$
where $\lacc \sigma_i, \, i\ge 1\racc$ is a deterministic sequence and $\lacc X_i, \, i\ge 1\racc$ an i.i.d. family. The situation is a bit simpler than for integrated random walks because the non-stationary increments of $\Sigma$ are independent, nevertheless the sequence $\lacc \sigma_i, \, i\ge 1\racc$ has also more generality. Setting $\mu$ for the law of $X_1$ and introducing
$$T_\sigma\; =\; \inf\lacc n\ge 1, \; \Sigma_n > 0\racc,$$
the following has been obtained in \cite{ABa} via strong approximation techniques:

\begin{THEO}[Aurzada-Baumgarten] Suppose that $\mu$ is centered and that its Laplace transform is finite in an open neighbourhood of zero. Suppose that $\lacc \sigma_i, \, i\ge 1\racc$ is increasing with $\sigma_n \asymp n^p$ for some $p > 0.$ Then
$$\pb [T_\sigma > n]\; =\; n^{-(p+1/2) +o(1)}.$$
\end{THEO}

The results of \cite{ABa} are more precise and allow other weight functions not necessarily increasing when $\mu$ is Gaussian. The general case reduces to the Gaussian one via strong approximation and drift transformations, as in \cite{AD}. A universal speed not depending on $\mu$ can also be obtained for weight functions growing faster than polynomials, like $e^{n^\gamma}$ for some $\gamma < 1/4,$ and the surviving probability is then exponentially small. However, it is also shown in \cite{ABa} that the speed does depend on $\mu$ for weight functions growing too fast, like $e^{n^\gamma}$ for some $\gamma \ge 1.$ It would be interesting to find $\sigma$'s critical growth rate for the universality of the speed. The following is also a natural question.

\begin{Conj} Suppose that $\mu$ is centered and has finite variance. Suppose that $\sigma_n \asymp n^p$ for some $p > 0.$ Then
$$\pb [T_\sigma > n]\; =\; n^{-(p+1/2) +o(1)}.$$
\end{Conj}

Let us conclude this paragraph with iterated L\'evy processes, which are processes of the type $\lacc X\circ \vert Y_t\vert, \; t\ge 0\racc$ with $X,Y$ two independent real L\'evy processes starting from zero. If $\vert Y\vert$ is a subordinator, then $X\circ \vert Y\vert$ is another L\'evy process which is called a subordinate L\'evy process, a notion introduced by Bochner in the context of harmonic analysis. Iterated L\'evy processes were introduced by Burdzy in the Brownian framework and can be viewed as a generalisation of subordinate L\'evy processes. They are known to have strong connections with PDE's of higher order, especially through their first passage times. Let us introduce
$$T_x\; =\; \inf\lacc t > 0, \; X\circ \vert Y_t\vert > 1\racc$$ 
for all $x \ge 0,$ and notice that as for integrated L\'evy processes the law of $T_x$ is difficult to study in general since $X\circ \vert Y\vert$ is non-Markov. Among other results, the following has been obtained in \cite{Ba}:

\begin{THEO}[Baumgarten] Suppose that the random variables $\vert X_1\vert^\a$ and $\vert Y_1\vert^\a$ have exponential moments for some $\a > 0$ and that $\EE [X_1] = 0.$ Then
$$\pb[T_1 > t]\; =\; t^{-\theta + o(1)}$$
with $\theta = 1/4$ if $\EE [Y_1] = 0$ and $\theta = 1/2$ if $\EE [Y_1] \neq 0.$ 
\end{THEO}
In particular, the persistence exponent of iterated Brownian motion is $1/4$. The strong dichotomy between the situations where $Y_1$ is centered and non-centered is not surprising since in the former case $\vert Y_t\vert $ grows roughly like $\sqrt{t}$ whereas in the latter case $\vert Y_t\vert $ grows like $t.$ As in many above statements, the question of replacing the exponential moment condition by the sole assumption of finite variance remains open.

\section{Some connections with physics}

\subsection{Regular points of inviscid Burgers equation with self-similar initial data} The statistical study of the one-dimensional Burgers equation
\begin{equation}
\label{Burgers}
\partial_t u \; +\; u\partial_x u \; =\; \nu \partial_{xx} u
\end{equation}
with viscosity $\nu > 0$ and an initial condition $u_0(x) := u(0,x) = X_x$ given by a self-similar stochastic process $\lacc X_x, x\in\rl\racc$ has been initiated in the papers \cite{AFS, Si2}. Though this equation is accorded to be an unrealistic physical model for turbulence, the competition between the irregularities of $X$ and the irregularities generated by (\ref{Burgers}) remains an interesting mathematical study. In the inviscid limit $\nu =0,$ the Hopf-Cole solution to (\ref{Burgers}) is given by
$$u(t,x)\; =\; \frac{x - a(t,x)}{t}$$
for every $t> 0, x\in \rl,$ where $a(t,x) = \max\{y\in\rl, \, {\dot C}_t(y)\le xt^{-1}\}$ and ${\dot C}$ is the right-derivative of the convex minorant of the function
$$F_t : y\; \mapsto\; \int_0^y (X_x + xt^{-1}) \d x.$$
This variational formula is obtained in considering the explicit solution to (\ref{Burgers}) which can be obtained for $\nu > 0$ and letting $\nu \to 0$ - see \cite{AFS, Be1, JW, Si2} for details. Notice that $a(t,x)$ is well-defined only if
\begin{equation}
\label{growth}
\vert x\vert^{-1}X_x\;\to\; 0 \;\;\;\mbox{a.s. when $x \to \pm\infty.$}
\end{equation}
The function $x\mapsto a(t,x)$ is right-continuous but not continuous in general, and the so-called Lagrangian regular points at time $t >0$ are defined as the set
$$\cL_t\; =\; \overline{\lacc a(t,x), \; x\in\rl\;\;\mbox{and}\;\; a(t,x-) =a(t,x)\racc}$$
which consists of the points where $F_t$ equals its convex minorant. In physical terms, the set $\cL_t$ describes the initial locations of the particles which have not been shocked up to time $t.$ It is easy to see that when $X$ is a self-similar process, the map $t\mapsto\cL_t$ has also some self-similarity which makes the a.s. Hausdorff dimension of $\cL_t$ independent of $t > 0.$ Setting $\cL =\cL_1$ and "Dim" for "Hausdorff dimension", the following is stated in \cite{AFS}:

\begin{Conj}[Aurell-Frisch-She] Suppose that $X$ is the fractional Brownian motion with Hurst parameter $H\in (0,1)$. Then {\em Dim} $\cL\, = H$ a.s. 
\end{Conj}

Notice that in the above, the fractional Brownian motion is defined over the whole $\rl,$ and coincides with the two-sided Brownian motion when $H = 1/2.$ This conjecture remains open in general, but the following has been shown:

\begin{THEO}[Handa-Sina\u {\i}, Bertoin] Suppose that $X$ is the two-sided Brownian motion. Then  {\em Dim} $\cL\, = 1/2$ a.s. 
\end{THEO}

This result was first stated in \cite{Si2}, although no strict proof is given therein for the lower bound Dim $\cL\, \ge 1/2$ a.s. A simple argument based on integration by parts and Frostman's lemma is presented in \cite{H}, which yields the general lower bound Dim $\cL\, \ge H$ a.s. when $X$ is the fractional Brownian motion with Hurst parameter $H.$ In \cite{Be1}, the exact computation of the Hausdorff dimension follows as a simple corollary to the more general result that $x\mapsto a(1,x)$ has stationary and independent increments with explicit Laplace transform. This result extends to L\'evy processes with no positive jumps satisfying the growth condition (\ref{growth}). In particular one has the

\begin{THEO}[Bertoin] Suppose that $X$ is a two-sided $\a$-stable spectrally negative L\'evy process with index $\a\in(1,2).$ Then {\em Dim} $\cL\, = 1/\a$ a.s. 
\end{THEO}

We now briefly describe the link between an upper bound for Dim $\cL$ and the computation of certain persistence exponents, in the self-similar framework. This is the original argument of \cite{Si2} for Brownian motion and it extends to fractional Brownian motion \cite{KM} or $\a$-stable L\'evy processes \cite{TS2}. Specifically, it can be shown by the Borel-Cantelli lemma and some elementary inequalities that if
\begin{equation}
\label{2sided}
\pb\lcr \int_0^y (X_x + x) \d x\, \ge \, -\delta^{1+H}, \; \forall y\in[-1,1]\rcr\; \le\; \delta^{1-K +o(1)}
\end{equation}
as $\delta\to 0,$ where $H$ is $X$'s self-similarity index and $K \in [0,1]$, then Dim $\cL\, \le K$ a.s. If $X=B^H$ is the fractional Brownian motion, then the drift appearing in (\ref{2sided}) can be removed by quasi-invariance and some analysis, and one sees by symmetry and self-similarity that the required estimate to get the upper bound in Conjecture 11 is
$$\pb\lcr \int_0^y B^H_x \d x\, \le \, 1, \; \forall y\in[-t,t]\rcr\; \le\; t^{1-H +o(1)}, \quad t\to +\infty.$$
This is a "two-sided" persistence problem, and the above estimate is formulated as a conjecture in \cite{KM, Mo2}, with an equality instead of the inequality. Notice that this latter problem is independent of Conjecture 9 since the increments of $B^H$ are correlated. Actually, even the sole existence of the persistence exponent for integrated double-sided fractional Brownian motion has not yet been established. \\
 
If $X$ is a two-sided $\a$-stable L\'evy process with $\a\in (1,2)$, some tedious analysis shows that the drift appearing in (\ref{2sided}) can also be removed (this would not be the case for $\a\in (0,1]$). By self-similarity and independence of the positive and negative increments of $Z$ the inequality (\ref{2sided}) amounts then to 
$$\pb[ {\hat T}_1 > t]\; \le\; t^{(1-K)/2 + o(1)}$$
at infinity, where ${\hat T}_1$ is the first-passage time at 1 of the integral of ${\hat Z} = -Z.$ In particular the validity of Conjecture 4 would lead to Dim $\cL\, \le \rho$ a.s. where $\rho = \pb[Z_1 > 0]$ is the positivity parameter of $Z$. In view of the above theorem, it is natural to state the
 
\begin{Conj} Suppose that $X$ is a two-sided $\a$-stable L\'evy process with $\a\in (1,2)$ and positivity parameter $\rho\in [1-1/\a, 1/\a]$. Then {\em Dim} $\cL\, = \rho$ a.s. 
\end{Conj}

This Hausdorff dimension depending on the positivity parameter and not on the self-similarity index is different from the value $1/\a$ which had been conjectured in \cite{JW} through multifractal analysis. The invalidity of the latter when $\a$ is close to 1 was proved in \cite{TS2} with the help of the above Sina\u {\i}'s approach:
 
\begin{THEO}[Simon] For every $c < 1$ there exists $\a_0 > 1$ such that for every $\a\in (1, \a_0)$ and every $\rho\in [1-1/\a, c\wedge1/\a],$ if $X$ is a two-sided $\a$-stable L\'evy process with positivity parameter $\rho$, then {\em Dim} $\cL\, < 1/\a$ a.s.  
\end{THEO} 

We conclude this paragraph in mentioning that the lower bound in Conjecture 12, which should be obtained by different and yet unknown arguments, would lead by comparison to the lower bound of Conjecture 4 in the case $\a > 1.$ 

\subsection{Positivity of random polynomials and diffusion equation} 

\label{sec:randpoly}

A classical question dating back to the beginning of probability theory is to understand the distribution of the roots of random polynomials. Consider 
$$P_n(X) \; = \;\sum_{i=0}^{n-1} \xi_i X^i$$
with large even degree where  $\lacc \xi_i, \; i\ge 0\racc$ is some i.i.d. sequence and $X$ the deterministic variable, and set $N_n$ for the number of its real roots. Among other results, the following was recently obtained in \cite{DPSZ}:

\begin{THEO}[Dembo-Poonen-Shao-Zeitouni] Suppose that $\xi_1$ is centered and has polynomial moments of all order. Then
\begin{equation} 
\label{eqn:DPSZ}
\pb \lcr N_{2n+1} = 0\rcr\; =\; n^{-4b +o(1)}
\end{equation}
where $b= - \lim_{t\to\infty}t^{-1}\log \pb\lcr Y_s < 0,\; \forall\, s\in[0,t] \rcr,$ with $Y$ the centered {\em GSP} with correlation $\EE [Y_0 Y_t] = 1/\cosh (t/2).$
\end{THEO}

In the above, the exact value of $b$ is unknown and numerical simulations suggest $4b =0.79 \pm 0.03$ - see \cite{DPSZ}. It is remarkable that this constant $b$ relates to $n$ times integrated Brownian motion. We saw indeed in Paragraph~\ref{sec:filp} that $b =  \lim_{n\to +\infty}\theta(n)$ where $\theta (n)$ is the persistence exponent of the process
$$t\; \mapsto\;\frac{1}{n!}\int_0^t (t-s)^n \, \d B_s, \quad t\ge 0.$$
The problem of computing $\theta(n)$ for $n > 1$ is believed to be very challenging. Numerical simulations \cite{BM} suggest $\theta(2) = 0.231\pm 0.01.$ \\

Let us give some insight on the proof of the above result. The hard part is to show (\ref{eqn:DPSZ}) when $\xi_1\sim \cN(0,1).$ The general case follows by strong approximation, whence the assumption made on the moments, but notice that it is also conjectured in \cite{DPSZ} that (\ref{eqn:DPSZ}) should hold under the sole condition that $\xi_1$ is centered and has finite variance. When $\xi_1\sim \cN(0,1)$, the process $x\to P_n(x),$ which is the so-called Kac's polynomial, is centered Gaussian with covariance $\EE [P_n(x)P_n(y)] = 1 +\cdots + (xy)^{n-1}.$ Its correlation function is given by
$$\lva\frac{(xy)^n - 1}{xy - 1}\rva\sqrt{\lva\frac{(x^2- 1)(y^2- 1)}{(x^{2n}- 1)(y^{2n}- 1)}\rva}$$
for every $x, y\neq \pm 1.$ This function is invariant under the transformations $(x,y)\mapsto (-x,-y)$ and $(x,y)\mapsto (1/x,1/y),$ and an involved argument based on Slepian's lemma shows that 
$$\pb \lcr N_{2n+1} = 0\rcr\; =\; ( \pb [N^{[0,1]}_{2n+1} = 0] )^{4+o(1)}$$
where $N^{[0,1]}_{n}$ is the number of roots of $P_n$ on $[0,1].$ The link between $N^{[0,1]}_{2n+1}$ for large $n$ and the zero-crossings of $Y$ is established after changing the variable $x=e^{-t}$ and isolating the contributions for small $t.$ The latter is the crucial step, since it follows from the singularities of the correlation function that the density of the real roots of $P_n$ around $\pm 1$ is very big for large $n.$ Notice that the link between $Y$ and $N_n$ is rather easily understood at the expectation level: a classical formula due to Kac - see \cite{DPSZ} for details - yields
$$\EE[N^{[0,1]}_{2n+1}]\; \sim\; \frac{1}{2\pi}\log n,$$
whereas the It\^o-Rice's formula shows that if $N_t$ is the number of zero-crossings of $Y$ on $[0,t],$ then $\EE[N_t]\sim t/2\pi.$ The problem of evaluating $\pb [N^{[0,1]}_{2n+1} = 0]$ is however much more intricate than the sole estimation of $\EE[N^{[0,1]}_{2n+1}],$ in analogy with what happens for the zero-crossings of Gaussian stationary processes. \\

It is also a remarkable fact that the above constant $b$ appears as the persistence exponent of another, apparently disconnected random evolution phenomenon, which is studied in \cite{M, MS, MS1}. Consider the heat equation on $\rl^d$
\begin{equation}
\label{Heat}
\frac{\partial u_d}{\partial t}\; =\; \Delta u_d
\end{equation}
with random initial condition $u_d(x,0) = {\dot W} (x)$ a $d$-dimensional white noise. Integrating along the heat kernel, it is easy to see that for every $x\in\rl^d$ the solution $t\mapsto u_d(x,t)$ to (\ref{Heat}) is a $(-d/4)$-self-similar centered Gaussian process with covariance function
$$\EE[u_d(t,x)u_d(s,x)]\; =\; \frac{1}{(\pi(t+s))^{d/2}}\cdot$$
In particular the law of $\lacc u_d(t,x), \; t > 0\racc$ does not depend on $x$, which is also clear from the white noise initial condition. The Lamperti transformation $t\mapsto (2\pi)^{1/4} e^{t/4} u_d(x, e^t)$ is a centered GSP with correlation function $1/(\cosh (t/2))^{d/2},$ and this GSP coincides with $Y$ for $d = 2.$ For every $d \ge 1,$ the standard superadditivity argument yields the existence of $\lbd_d > 0$ such that
$$\pb[u_d(s,x) < 0, \; \forall s\in [1,t]]\; =\; t^{-\lbd_d + o(1)}$$
and one has $b = \lbd_2.$ In \cite{M}, an empirical approach using independent interval approximation is described, proposing $\lbd_d$ as the first zero on the negative axis of the function
$$x\; \mapsto\; 1 \; +\; \sqrt{\frac{2}{d}}\lpa \pi x - 2 x^2\int_0^\infty e^{-xt} \sin^{-1}[1/(\cosh (t/2))^{d/2}] dt\rpa,$$
which yields the numerical values $\lbd_1 = 0.1207, \lbd_2 = 0.1862$ and $\lbd_3 = 0.2358.$ In the above formula, the function $d \mapsto \lbd_d$ is increasing, which is somehow in heuristic accordance with the fact that in the first-order expansion
$$1/(\cosh (t/2))^{d/2}\; =\; 1\; -\; \frac{d t^2}{16} \; +\; o(t^2),$$
the coefficient $d \mapsto d/16$ also increases. In \cite{MS} it is argued that $\lbd_d\sim c\sqrt{d}$ at infinity, for some constant $c >0.$ The paper \cite{MS} also establishes for every $d \ge 1$ a general connection between the survival analysis of the equation (\ref{Heat}) and the positivity of a family of random polynomials defined as
$$P^d_n(X) \; = \;\xi_0\; +\; \sum_{i=1}^{n-1} i^{(d-2)/4}\xi_i X^i$$
where $\lacc \xi_i, \; i\ge 0\racc$ is a i.i.d. sequence of $\cN(0,1)$ random variables and $X$ is the deterministic variable. Setting $N^d_n$ for the number of its real roots, it is argued in \cite{MS} that
\begin{equation}
\label{GKac}
\pb \lcr N^d_{2n+1} = 0\rcr\; =\; n^{-2(\lbd_2 +\lbd_d) +o(1)}.
\end{equation}
The method relies on an analysis of the mean density $\rho_n^d$ of the real roots, which is defined as
$$\EE[N^d_n [a,b]]\; =\; \int_a^b \rho_n(x) dx$$
for every $a<b$, where $N^d_n[a,b]$ is the number of real roots of $P^d_n(X)$ on $[a,b]$. 
It is proved in \cite{MS} that the limit mean density $\rho_\infty^d$ has a shape which is independent of $d$ on $(-1, 1)^c$ (whence the contribution $2\lbd_2$ in the exponent), and depends on $d$ on $(-1,1)$ (whence the contribution $2\lbd_d$). The rate of explosion at $\pm 1$ is also explicitly evaluated, which makes it possible to obtain the precise link between $N^d_n(-1,1)$ and the process $u_d$, resp. between $N^d_n(-1,1)^c$ and the process $u_2$. We confess not having checked the whole argument of \cite{MS} in detail, which seems to deserve a more complete mathematical explanation as it is done in \cite{DPSZ} for the case $d = 2$. 

\subsection{Wetting models with Laplacian interactions} Let $f$ be a bounded and everywhere positive probability density over $\rl$, centered and having finite variance. Introduce the Hamiltonian $\cH_{[a,b]}(\varphi),$ defined for $a, b \in \ZZ$ with $b-a \ge 2$ and for $\varphi : \{a, \ldots, b\}\to\rl$ by
$$\cH_{[a,b]}(\varphi)\; =\; \sum_{n = a+1}^{b-1} V(\Delta \varphi_n)$$
where $V = -\log (f)$ is the potential and 
$$\Delta \varphi_n\; =\; (\varphi_{n+1} - \varphi_n) - (\varphi_n - \varphi_{n-1})\; = \; \varphi_{n+1} - 2\varphi_n + \varphi_{n-1}$$
is the discrete Laplacian on $\ZZ.$ The free pinning model with Laplacian interaction is the probability measure on $\rl^{N-1}$ defined by
$$\pb_{0, N}^p(\d \varphi_1, \ldots, \d \varphi_{N-1})\; =\; \frac{\exp (-\cH_{[-1,N+1]}(\varphi))}{\cZ^p_{0,N}}\d \varphi_1 \ldots \d \varphi_{N-1}$$
where $\cZ^p_{0,N}$ is the normalization constant which is called the partition function, and where the boundary conditions are given by $\varphi_{-1} = \varphi_{0} = \varphi_{N} = \varphi_{N+1} = 0.$ This probability measure modelises a certain $(1+1)-$dimensional field (viz. a linear chain $\{(n, \varphi_n), \, n = 0\ldots N\}$) with zero boundary conditions and whose interacting structure is described by the discrete Laplacian and the potential $V$. This chain can be viewed as an example of a discrete random polymer in $(1+1)-$dimension.

The free pinning model with gradient interaction, where $\Delta$ is replaced by the discrete gradient $\nabla \varphi_n = \varphi_{n+1} - \varphi_n$ and where the boundary conditions are $\varphi_{0} = \varphi_{N} = 0,$ has been well studied in the literature and has a natural interpretation in terms of random bridges with increment density given by $f$. The model with Laplacian interaction has exactly the same interpretation in terms of integrated random bridges. Specifically, one can easily show - see Section 2 in \cite{CD} for details - that $\pb_{0, N}^p$ is, with the notations of Section \ref{sec:irw}, the law of an integrated random walk $\{A_n= S_1 +\cdots +S_n, \, n = 1\ldots N-1\}$ conditioned on $A_N = A_{N+1} = 0,$ with increment density given by $f$. The partition function $\cZ^p_{0,N}$ is then the value at $(0,0)$ of the density of $(S_{N+1}, A_{N+1}).$ Notice that both above free pinning models have natural counterparts in continuous time in the context of semiflexible polymers. The gradient interacting case corresponds to directed polymers, whereas the Laplacian interacting case corresponds to polymers with non-zero bending energy - see \cite{Bu} for details.

The connection with persistence of integrated random bridges is made in considering the corresponding wetting model with Laplacian interactions, which is the probability measure on $\rl^{N-1}$ defined by
\begin{eqnarray*}
\pb_{0, N}^w(\d \varphi_1, \ldots, \d \varphi_{N-1})& =& \pb_{0, N}^p(\d \varphi_1, \ldots, \d \varphi_{N-1}\;\vert\; \varphi_1\ge 0, \ldots, \varphi_{N-1}\ge 0)\\
& =& \frac{\exp (-\cH_{[-1,N+1]}(\varphi))}{\cZ^w_{0,N}}\Un_{\{\varphi_1\ge 0, \ldots, \varphi_{N-1}\ge 0\}}\d \varphi_1 \ldots \d \varphi_{N-1}
\end{eqnarray*}
where $\cZ^w_{0,N}$ is the normalization constant, and with the same boundary conditions $\varphi_{-1} = \varphi_{0} = \varphi_{N} = \varphi_{N+1} = 0.$ In this model, the discrete random polymer is in the presence of a one-dimensional hard wall at zero which forces it to stay non negative. From the very definition, one sees that $\pb_{0, N}^W$ is the law of the above integrated random walk $\{A_n, \, n = 1\ldots N-1\}$ whose increment has the density $f$, and conditioned on $\Omega_{N-1}^+\cap\{A_N = A_{N+1} = 0\}$ with the notation $\Omega_{N-1}^+ = \{A_1\ge 0, \ldots, A_{N-1}\ge 0\}.$ It is then easily shown that the partition function is given by
$$\cZ^w_{0,N}\; =\; \pb[\Omega_{N-1}^+\, \vert\, A_N = A_{N+1} = 0] f_N(0,0)$$
where $f_N(0,0)$ is the value at $(0,0)$ of the density of $(S_{N+1}, A_{N+1})$. As a consequence of a local limit theorem - see Section 2 in \cite{CD} for details - it can be shown that $f_N(0,0)\sim c N^{-2}$ at infinity for some explicit constant $c > 0.$ Hence the behaviour of $\cZ^w_{0,N}$ for $N$ large, which has some importance in physics, is specified by the persistence probability 
$$\pb[\Omega_{N-1}^+\, \vert\, A_N = A_{N+1} = 0].$$ 
The latter quantity has also some independent interest as a question about entropic repulsion - see all the references listed in \cite{CD} for more on this subject, and the following is stated in \cite{CD}:

\begin{Conj}[Caravenna-Deuschel] With the above notations, one has
$$\pb[\Omega_{N-1}^+\, \vert\, A_N = A_{N+1} = 0]\; \asymp\; N^{-1/2}$$
for every centered increment law $\mu$ having finite variance. 
\end{Conj}

This conjecture is related with Conjecture 1 since the event $\{\Omega_{N-1}^+\, \vert\, A_N = A_{N+1} = 0\}$ can be decomposed into 
$\{A_1\ge 0, \ldots, A_{N/2}\ge 0\, \vert\, A_N = A_{N+1} = 0\}\cap\{A_{N/2+1}\ge 0, \ldots, A_{N-1}\ge 0\, \vert\, A_N = A_{N+1} = 0\},$
the intersection of two roughly independent events with roughly the same probability $\pb[A_1\ge 0, \ldots, A_{N/2}\ge 0],$ a quantity which should behave like $N^{-1/4}.$ Notice that in the context of semiflexible polymers, a continuous couterpart of $\Omega_{N-1}^+$ in the case when $\mu$ is Gaussian was investigated (without conditioning) in \cite{Bu}, where the estimate (\ref{GIW}) is proved. In \cite{CD}, the following weak bounds are obtained 
$$\frac{c}{N^{c_-}}\; \le\; \pb[\Omega_{N-1}^+\, \vert\, A_N = A_{N+1} = 0]\; \le\; \frac{C}{(\log N)^{c_+}}$$
for some constants $c, C, c_- >0$ and $c_+ > 1.$ The lower bound entails that the free energy vanishes:
$$\lim_{N\to +\infty}\frac{1}{N}\log\cZ^w_{0,N}\; =\; 0,$$
whereas the fact that $c_+ > 1$ in the upper bound is crucial to show that the phase transition of the wetting model with reward, which is simply the value of the positive parameter $\varepsilon$ after which the free energy of the probability measure
$$\frac{\exp (-\cH_{[-1,N+1]}(\varphi))}{\cZ^w_{\varepsilon,N}}\prod_{n=1}^{N-1}(\varepsilon \delta_0(\d \varphi_n) + \Un_{\{\varphi_n\ge 0\}}\d \varphi_n) $$
becomes positive, is of first order. In \cite{CD}, to which we again refer for more details, it is mentioned that Conjecture 13 would yield some further path results for the wetting model with reward at criticality.

\subsection{Other physical applications}

\subsubsection{Spatial persistence for fluctuating interfaces} A fluctuating interface is a function $h : \rl^+\times \rl^d\to \rl$ evolving in time, with dynamics governed by a certain random equation. The problem of spatial persistence concerns the probability $p(l)$ that such a fluctuating interface stays above its initial value over a large distance $l$ from a given point in space. One expects a behaviour like $p(l) = l^{-\theta +o(1)}$ for a positive number $\theta$  independent of the direction, which is called the spatial persistence of the interface. In \cite{BM}, this question is adressed for the Gaussian interface $h(t,x)$ solution to the equation
\begin{equation}
\label{Lange}
\frac{\partial h}{\partial t}\; =\; - (-\Delta)^{z/2}h \; +\; \xi
\end{equation}
where $\Delta$ is the $d-$dimensional Laplacian, $\xi$ a space-time white noise with zero mean, and $z > d$ some fractional parameter, and it is shown with heuristic arguments based on Fourier inversion that if $\rho = (z-d +1)/2$, then the fractional derivative in any direction $x_1$
$$\frac{\partial^\rho h}{\partial x_1^\rho}\; =\; \frac{\partial}{\partial x_1^{[\rho]+1}}\lpa \frac{1}{\Gamma(2-\a)}\int_0^{x_1}h(t, y, x_2, \ldots, x_d)(x_1-y)^{1-\a} \d y\rpa,$$
where we have decomposed $\rho = [\rho] + \a$ into integer and fractional parts, is a one-dimensional white noise. This relates the spatial persistence probability $p(l)$ to the persistence probability of the Riemann-Liouville process
$$A^\rho_t \; =\; \frac{1}{\Gamma(\rho)}\int_0^{t}(t-s)^{\rho-1} \d B_s.$$
In \cite{BM} two different regimes are considered. The coarsening one, where the reference point is fixed, yields a spatial persistence exponent $\theta = \theta (\rho)$ with the notations of Paragraph \ref{sec:filp}. The stationary one, where the reference point is sampled uniformly from the ensemble of steady state configurations, yields from (\ref{Incr}) a spatial persistence exponent $\theta = (1- \rho)_+.$ In the coarsening regime, this entails that the zero crossings of Gaussian interfaces governed by (\ref{Lange}) undergo a morphological transition at $z = d+2,$ because then $\rho =3/2.$

\subsubsection{Clustering of sticky particles at critical time}  In this last paragraph we consider a random walk $\{S_n, \, n\ge 1\}$ with positive increments having expectation $\EE[S_1] = 1.$ If $n$ particles are fixed at the respective positions $i^{-1}S_i, i = 1\ldots n$ with zero initial speed and then move according to the laws of  gravitational attraction, these particles end up in sticking together with conservation of mass and momentum, forming new particles called clusters. One is then interested in the number of clusters $K_n(t)\in [1, \ldots, n]$ viz. the total number of particles present at time $t\geq 0$. This is a so-called sticky particle model, which is for $n$ large connected to the inviscid Burgers equation with random initial data (coupled with some scalar transport equation, see \cite{BrG} for details). This is also an aggregation model having connections with astrophysics, and we refer to the introduction of \cite{Vy0} for a clarification of these relations and the complete dissipation of all possible misunderstanding.

The normalization $\EE[S_1] = 1$ entails that $T\to 1$ in probability, where $T$ stands for the random terminal time where all particles have aggregated in a single cluster. A more precise result is obtained in \cite{Vy0} in the case when $S_1$ has uniform or standard Poissonian distribution, namely that the random function
$$\frac{K_n(t) - n (1-t^2)}{\sqrt{n}}$$
converges in law to some Gaussian process on the Skorokhod space $\cD[0,1-\eps]$ for any $\eps>0.$ In particular, the above quantity converges to some Gaussian law at each fixed time $t<1.$ The situation is however different at the critical time $t=1,$ at least when $S_1$ has a standard Poissonian distribution. In this case it can be proved that $K_n(1)/\sqrt{n}$ does not converge to zero (the only non-negative Gaussian distribution) as could be expected, and this fact is actually a consequence - see \cite{Vy0} for details - of the estimate 
$$\pb\lcr \min_{i=1,\ldots, n} \sum_{j=1}^i (\Gamma_j-j) \geq 0\rcr \; \asymp\; n^{-1/4}$$
where $\{\Gamma_n, \, n\ge 1\}$ is a random walk with exponential increments - this latter estimate follows from the main result of \cite{DG}.

\bigskip

\noindent
{\bf Acknowledgements.} This work was supported by the DFG Emmy Noether research program and the grant ANR-09-BLAN-0084-01. The authors would like to thank R.~A.~Doney, S.~N.~Majumdar, and A.~A.~Novikov for some useful comments.

\end{document}